\documentclass[12pt,oneside,reqno]{amsart}
\usepackage{mathrsfs}
\usepackage{graphics}
\usepackage{amssymb}
\usepackage{enumerate}
\newcommand{\dif}{\mathrm{d}}

\newcommand{\be}{\begin{eqnarray}}
\newcommand{\ee}{\end{eqnarray}}
\newcommand{\ce}{\begin{eqnarray*}}
\newcommand{\de}{\end{eqnarray*}}
\newtheorem{theorem}{Theorem}[section]
\newtheorem{lemma}[theorem]{Lemma}
\newtheorem{remark}[theorem]{Remark}
\newtheorem{definition}[theorem]{Definition}
\newtheorem{proposition}[theorem]{Proposition}
\newtheorem{Example}[theorem]{Example}
\newtheorem{corollary}[theorem]{Corollary}
\newtheorem{condition}[theorem]{Condition}
\def\e{\varepsilon}

\def\a{\alpha}

\def\d{\delta}

\def\g{\gamma}

\def\[{{\Big[}}
\def\]{{\Big]}}
\def\<{{\langle}}
\def\>{{\rangle}}
\def\({{\Big(}}
\def\){{\Big)}}

\def\no{\nonumber}
\def\bt{\begin{theorem}}
\def\et{\end{theorem}}
\def\bl{\begin{lemma}}
\def\el{\end{lemma}}
\def\br{\begin{remark}}
\def\er{\end{remark}}
\def\bx{\begin{Example}}
\def\ex{\end{Example}}
\def\bd{\begin{definition}}
\def\ed{\end{definition}}
\def\bp{\begin{proposition}}
\def\ep{\end{proposition}}
\def\bc{\begin{corollary}}
\def\ec{\end{corollary}}
\def\bco{\begin{condition}}
\def\eco{\end{condition}}

\def\cD{{\mathcal D}}

\def\cG{{\mathcal G}}

\def\cP{{\mathcal P}}

\def\mE{{\mathbb E}}

\def\mH{{\mathbb H}}

\def\mN{{\mathbb N}}

\def\mP{{\mathbb P}}

\def\mR{{\mathbb R}}
\def\mS{{\mathbb S}}

\def\mW{{\mathbb W}}

\def\sA{{\mathscr A}}
\def\sB{{\mathscr B}}

\def\sF{{\mathscr F}}

\def\sL{{\mathscr L}}

\def\sV{{\mathscr V}}

\def\geq{\geqslant}
\def\leq{\leqslant}
\def\epsilon{\varepsilon}
\begin{document}

\allowdisplaybreaks

\title{Deviation estimates  for multivalued McKean-Vlasov stochastic differential equations}

\author{Kun Fang and Huijie Qiao$^*$}

\dedicatory{School of Mathematics,
Southeast University\\
Nanjing, Jiangsu 211189, China}



\thanks{This work was partly supported by NSF of China (No. 12071071).}

\thanks{*Corresponding author: hjqiaogean@seu.edu.cn}


\date{}

\begin{abstract}
The work concerns deviation estimates for multivalued McKean-Vlasov stochastic differential equations. First of all, we prove the large deviation principle for them by the weak convergence approach. Then the central limit theorem for them is shown with the help of a formula for $L$-derivatives. Finally, we establish the moderate deviation principle for them.
\end{abstract}

\maketitle \rm

\section{Introduction}
Given a filtered probability space $(\Omega,\mathscr{F},\{\mathscr{F}_t\}_{t\geq 0},\mP)$ and a $m$-dimensional standard Brownian motion $W_{\cdot}=(W_{\cdot}^1,W_{\cdot}^2,\cdots, W_{\cdot}^m)$ defined on it. Consider the following multivalued McKean-Vlasov stochastic differential equation (SDE for short) on $\mR^d$:
\be\left\{\begin{array}{l}
\dif X_t\in \ -A(X_t)\dif t+ \ b(X_t,\sL_{X_t})\dif t+\sigma(X_t,\sL_{X_t})\dif W_t,\\
X_0=\xi\in\overline{\cD(A)}, \sL_{X_t}=\mP_{X_t}= \mbox{the probability  distribution of}~X_t,
\end{array}
\label{eq1}
\right.
\ee
where $\xi$ is nonrandom, $A:\mR^d \mapsto 2^{\mR^d} $ is a maximal monotone operator, $Int(\cD(A))\ne\emptyset$ and the coefficients $b:\mR^d\times\cP_2(\mR^d)\mapsto{\mR^d}, \,\,\sigma:\mR^d\times\cP_2(\mR^d)\mapsto{\mR^d}\times{\mR^m}$ are Borel measurable.

If the coefficients $b, \sigma$ don't depend on distributions of solution processes, Eq.(\ref{eq1}) is called a multivalued SDE. The type of equations was firstly introduced by C\'epa in \cite{cepa} and \cite{cepaa}, and then has received increasing attentions from researchers in recent years (c.f. \cite{Ren, rwz, RWZ1, RXZ, zh, ZXCH}). Let us mention some works related with ours. In \cite{RXZ}, Ren, Xu and Zhang  proved the large deviation principle for multivalued SDEs. Later, Ren, Wu and Zhang \cite{rwz} presented a general large deviation principle for them. Recently, Zhang \cite{zh} established the moderate deviation principle of them.

If the coefficients $b, \sigma$ depend on distributions of solution processes, Eq.(\ref{eq1}) is called a multivalued McKean-Vlasov SDE. The type of equations is the generalization of McKean-Vlasov SDEs and contains McKean-Vlasov stochastic variational inequalities where the maximal monotone operator $A$ is the subdifferential operator of some convex function. Although there are many results about McKean-Vlasov SDEs (c.f. \cite{DQ, DQ2, rst, Kac, lq, lszz, Re, Suo}), only a few results on multivalued McKean-Vlasov SDEs appear. Let us review some results. In \cite{CHI}, Chi proved the existence and uniqueness for the solutions of a type of multivalued McKean-Vlasov SDEs where the coefficients $b, \sigma$ depend on distributions through integrations. And Ren and Wang \cite{RW} showed the well-posedness and the uniform large deviation principle for mean-field stochastic variational inequalities where the coefficients $b, \sigma$ contain mathematical expectations. Very recently, Gong and Qiao \cite{G} established  the well-posedness and stability for Eq.(\ref{eq1}) under non-Lipschitz conditions.

In the paper, we follow up the line in \cite{G} and study the asymptotic behavior of the strong solution for Eq.(\ref{eq1}) in different deviation scales. Concretely speaking, for any $\e>0$, consider the following  multivalued McKean-Vlasov SDE:
\be\left\{\begin{array}{l}
\dif X^{\epsilon}_{t}\in -A(X^{\epsilon}_{t})\dif t+b(X^{\epsilon}_{t},\sL_{X^{\epsilon}_{t}})\dif t+\sqrt{\epsilon}\sigma(X^{\epsilon}_{t},\sL_{X^{\epsilon}_{t}})\dif W_t, \\
 X^{\epsilon}_{0}=\xi\in\overline{\cD(A)}.
\end{array}
\label{meq1}
\right.
\ee
Assume that $(X^{\epsilon}, K^{\epsilon})$ is a strong solution of Eq.(\ref{meq1}). Then we investigate the deviations of $X^{\epsilon}$ from $X^{0}$  by studying the asymptotic behavior of the trajectory
$$
\frac{X^{\epsilon}- X^{0}}{a(\epsilon)},
$$
where $a:\mathbb{R}^{+}\mapsto (0,1)$ and $(X^{0}, K^0)$ satisfies the following multivalued differential equation
\be\left\{\begin{array}{l}
\dif X_t^{0}\in \ -A(X_t^{0})\dif t+  b(X_t^{0},\delta_{X_t^{0}})\dif t,\\
X^0_0=\xi.
\end{array}
\label{eq3}
\right.
\ee
Our contribution is as follows:
\begin{enumerate}[$\bullet$]
	\item we provide the large deviation estimate for Eq.(\ref{eq1}) in the case of $a(\epsilon)\equiv1$.
\end{enumerate}
\begin{enumerate}[$\bullet$]
	\item we show the central limit theorem for Eq.(\ref{eq1}) in the case of $a(\epsilon)=\sqrt{\epsilon}$. That is, $\frac{X^{\epsilon}- X^{0}}{a(\epsilon)}$ converges to a stochastic process in a certain sense as $\epsilon\rightarrow0$. 
\end{enumerate}
\begin{enumerate}[$\bullet$]
	\item we prove the moderate deviation principle for Eq.(\ref{eq1}) in the case of $a(\epsilon)$ satisfying
    \ce
	a(\epsilon)\rightarrow0,\quad \frac{\epsilon}{a^{2}(\epsilon)}\rightarrow0\quad as\quad \epsilon\rightarrow0.
	\de
\end{enumerate}

It is worthwhile to mentioning our results and methods. Firstly, since our equation is more general than one in \cite{RW}, our result can cover \cite[Theorem 4.3]{RW}. Secondly, if the maximal monotone operator $A$ is zero, our equation becomes 
a McKean-Vlasov SDE. In \cite{Suo}, Suo and Yuan established the central limit theorem and the moderate deviation principle for McKean-Vlasov SDEs. Therefore, our result is more general. Thirdly, the traditional method of large and moderate deviation principles requires exponential tightness estimate and other exponential probability estimations (c.f. \cite{rst}). However, this method is particularly cumbersome for multivalued McKean-Vlasov SDEs. Hence, we use the weak convergence method to prove large and moderate deviation principles (c.f. \cite{de, lwyz, V, m}). Fourthly, in order to obtain the moderate deviation principle for Eq.(\ref{eq1}), because the coefficients $b, \sigma$ contain the distribution of the solution process, it is difficult directly to prove that $\frac{X^{\epsilon}- X^{0}}{a(\epsilon)}$ satisfies the large deviation 
principle (c.f.\cite{zh}). To overcome the difficulty, we show the moderate deviation principle for Eq.(\ref{eq1}) through the exponential equivalence.

The rest of this paper is organized as follows. In Section \ref{fram}, we recall some notions and some lemmas. Then a  general criterion of large deviation principles is given and the uniform large deviation principle for Eq.(\ref{eq1}) is proved. In Section \ref{clt}, we establish the central limit theorem for Eq.(\ref{eq1}). Finally, in Section \ref{mdp}, the moderate deviation principle for Eq.(\ref{eq1}) is presented.

The following convention will be used throughout the paper: $C$ with or without indices will denote different positive constants whose values may change from one place to another.

\section{Preliminary}\label{fram}

In the section, we introduce notations and concepts and recall some results used in the sequel.

\subsection{Notations}\label{nota}

In the subsection, we introduce some notations.

Let $\mid\cdot\mid$ and $\|\cdot\|$ be norms of vectors and matrices, respectively. Furthermore, let $\langle\cdot$ , $\cdot\rangle$ denote the scalar product in $\mR^d$. Let $B^*$ denote the transpose of the matrix $B$.

Let $C(\mR^d)$ be the collection of continuous functions on $\mR^d$ and $C^2(\mR^d)$ be the space of continuous functions on $\mR^d$ which have continuous partial derivatives of order up to $2$. 

Let $\sB(\mR^d)$ be the Borel $\sigma$-algebra on $\mR^d$ and $\cP({\mR^d})$ be the space of all probability measures defined on $\sB(\mR^d)$ carrying the usual topology of the weak convergence. Let $\cP_{2}(\mR^d)$ be the set of probability measures on $\sB(\mR^d)$ with finite second order moments. That is,
$$
\cP_2\left( \mathbb{R}^d \right) :=\left\{ \mu \in \cP\left( \mathbb{R}^d \right): \|\mu\| _{2}^{2}:=\int_{\mathbb{R}^d}{\left| x \right|^2\mu \left( \dif x \right) <\infty} \right\}.
$$
As we can see, $\cP_2(\mR^d)$ is a Polish space under the Wasserstein distance
$$
\mathbb{W}_2(\mu,\nu):= \inf\limits_{\pi\in\Psi(\mu,\nu)}\left(\int_{\mathbb{R}^d\times\mathbb{R}^d}|x-y|^{2}\pi(\dif x,\dif y)\right)^{\frac{1}{2}}, \quad \mu , \nu\in \cP_2(\mR^d),
$$
where $\Psi(\mu,\nu)$ is the set of couplings for $\mu$ and $\nu$.

\subsection{Maximal monotone operators}

In the subsection, we introduce maximal monotone operators.

Fix a multivalued operator $A: \mR^d\mapsto 2^{\mR^d}$, where $2^{\mR^d}$ stands for all the subsets of $\mR^d$, and set
$$
\cD(A):= \left\{x\in \mR^d: A(x) \ne \emptyset\right\}
$$
and
$$
Gr(A):= \left\{(x,y)\in \mR^{2d}:x \in \cD(A), ~ y\in A(x)\right\}.
$$
Then we say that $A$ is monotone if $\langle x_1 - x_2, y_1 - y_2 \rangle \geq 0$ for any $(x_1,y_1), (x_2,y_2) \in Gr(A)$, and $A$ is maximal monotone if
$$
(x_1,y_1) \in Gr(A) \iff \langle x_1-x_2, y_1 -y_2 \rangle \geq 0, \forall (x_2,y_2) \in Gr(A).
$$

Given $T>0$. Let $\sV_{0}$ be the set of all continuous functions $K: [0,T]\mapsto\mR^{d}$ with finite variations and $K_{0} = 0$. For $K\in\sV_0$ and $s\in [0,T]$, we shall use $|K|_{0}^{s}$ to denote the variation of $K$ on [0,s]. Set
\ce
&&\sA:=\Big\{(X,K): X\in C([0,T],\overline{\cD(A)}), K \in \sV_0, \\
&&\qquad\qquad\quad~\mbox{and}~\langle X_{t}-x, \dif K_{t}-y\dif t\rangle \geq 0 ~\mbox{for any}~ (x,y)\in Gr(A)\Big\}.
\de
Then about $\sA$ we recall the following results (cf.\cite{cepaa,ZXCH}).

\bl\label{equi}
For $X\in C([0,T],\overline{\cD(A)})$ and $K\in \sV_{0}$, the following statements are equivalent:
\begin{enumerate}[(i)]
\item $(X,K)\in \sA$.
\item For any $(x,y)\in C([0,T],\mR^d)$ with $(x_t,y_t)\in Gr(A)$, it holds that
$$
\left\langle X_t-x_t, \dif K_t-y_t\dif t\right\rangle \geq0.
$$
\item For any $(X^{'},K^{'})\in \sA$, it holds that
$$
\left\langle X_t-X_t^{'},\dif K_t-\dif K_t^{'}\right\rangle \geq0.
$$
\end{enumerate}
\el

\bl\label{inteineq}
Assume that $Int(\cD(A))\ne\emptyset$, where $Int(\cD(A))$ denotes the interior of the set $\cD(A)$. For any $a\in Int(\cD(A))$, there exists constants $\gamma_1 >0$, and $\gamma_{2},\gamma_{3}\geq0$ such that  for any $(X,K)\in \sA$ and $0\leq s<t\leq T$,
$$
\int_s^t{\left< X_r-a, \dif K_r \right>}\geq \gamma_1\left| K \right|_{s}^{t}-\gamma _2\int_s^t{\left| X_r-a\right|}\dif r-\gamma_3\left( t-s \right) .
$$
\el

\subsection{Multivalued McKean-Vlasov SDEs}

In the subsection, we introduce multivalued McKean-Vlasov SDEs.

Consider Eq.(\ref{eq1}), i.e.
\ce\left\{\begin{array}{l}
\dif X_t\in \ -A(X_t)\dif t+ \ b(X_t,\sL_{X_t})\dif t+\sigma(X_t,\sL_{X_t})\dif W_t,\\
X_0=\xi\in\overline{\cD(A)}, \sL_{X_t}=\mP_{X_t}= \mbox{the probability  distribution of}~X_t.
\end{array}
\right.
\de
A strong solution of Eq.(\ref{eq1}) means that there exists a pair of adapted processes $(X,K)$ on $(\Omega, \mathscr{F}, \{\mathscr{F}_t\}_{t\in[0,T]}, \mP)$ such that

(i) $\mP(X_0=\xi)=1$,

(ii) $X_t\in{\mathscr{F}_t^W}$, where $\{\mathscr{F}_t^W\}_{t\in[0,T]}$ stands for the $\sigma$-field filtration generated by $W$,

(iii) $(X_{\cdot}(\omega),K_{\cdot}(\omega))\in \sA$ a.s. $\mP$,

(iv) it holds that
\ce
\mP\left\{\int_0^T(\mid{b(X_s,\sL_{X_s})}\mid+\|\sigma(X_s,\sL_{X_s})\|^2)\dif s<+\infty\right\}=1,
\de
and
\ce
X_t=\xi-K_{t}+\int_0^tb(X_s,\sL_{X_s})\dif s+\int_0^t\sigma(X_s,\sL_{X_s})\dif W_s, \quad 0\leq{t}\leq{T}.
\de

\subsection{$L$-derivatives for the functions on $\cP_2(\mR^d)$}

In the subsection, we introduce $L$-derivatives for the functions on $\cP_2(\mR^d)$.

For any $\mu\in \cP_2(\mR^d)$, set
\ce
&&T_{\mu,2}:=L^{2}(\mathbb{R}^d\mapsto\mathbb{R}^d;\mu)\\
&&\qquad :=\Big\{\phi:\mathbb{R}^d\mapsto\mathbb{R}^d; \phi ~\mbox{is measurable with}~\mu(|\phi|^{2}):=\int_{\mR^d}|\phi(x)|^{2}\mu(\dif x)<\infty \Big\},\\
&&\|\phi\|_{T_{\mu,2}}^{2}:=\int_{\mathbb{R}^d}|\phi|^{2}\mu(\dif x), ~\mbox{for}~\phi\in T_{\mu,2}.
\de

\bd\label{compleve}
Let $f:\cP_2(\mR^d)\mapsto\mR$ be a continuous function, and $I$ be the identity map on $\mR^d$. 

(i) If for any $\mu\in\cP_2(\mR^d)$ 
$$
T_{\mu,2}\ni\phi\mapsto D_{\phi}^{L}f(\mu):=\lim\limits_{\epsilon\rightarrow 0}\frac{f(\mu\circ(I+\epsilon\phi)^{-1})-f(\mu)}{\epsilon}\in\mR
$$
is a well-defined bounded linear functional, we call $f$ intrinsically differentiable at $\mu$, and the intrinsic derivative of $f$ at $\mu$ is $D_{\phi}^{L}f(\mu)$.

(ii) If for any $\mu\in\cP_2(\mR^d)$ 
$$
\lim\limits_{\|\phi\|_{T_{\mu,2}}\rightarrow 0}\frac{f(\mu\circ(I+\epsilon\phi)^{-1})-f(\mu)-D_{\phi}^{L}f(\mu)}{\epsilon}=0,
$$
we call $f$ $L$-differentiable at $\mu$ and the $L$-derivative (i.e. Lions derivative) of $f$ at $\mu$ is denoted as $D^{L}f(\mu)$.
\ed

By the Riesz representation theorem, we know that
$$
\left\langle D^{L}f(\mu),\phi \right\rangle_{T_{\mu,2}}:=\int_{\mR^d}\left\langle D^{L}f(\mu)(x),\phi(x) \right\rangle\mu(\dif x)=D_{\phi}^{L}f(\mu), \quad \phi\in T_{\mu,2}.
$$

\bd
$f\in C^{1}(\cP_2(\mR^d))$ means that $f$ is $L$-differentiable at any point $\mu\in\cP_2(\mR^d)$, and the
L-derivative $D^{L}f(\mu)(x)$ has a version jointly continuous in $(\mu,x)\in\cP_2(\mR^d)\times\mR^d$. If
moreover, $D^{L}f(\mu)(x)$ is bounded, we denote $f\in C_{b}^{1}(\cP_2(\mR^d))$.
\ed

For a vector-valued function $f=(f_{i})$, or a matrix-valued function $f=(f_{ij})$ with $L$-differentiable components, we write
$$
D_{\phi}^{L}f(\mu)=(D_{\phi}^{L}f_{i}(\mu)),\quad or \quad D_{\phi}^{L}f(\mu)=(D_{\phi}^{L}f_{ij}(\mu)),\quad \mu\in\cP_2(\mR^d).
$$

\section{The large deviation principle for multivalued McKean-Vlasov SDEs}\label{lde}

In this section, we study the large deviation principle for multivalued McKean-Vlasov SDEs.

\subsection{A general criterion of large deviation principles}

In the subsection, we introduce a general criterion of large deviation principles.

The theory of small-noise large deviations concerns with the asymptotic behavior of solutions of multivalued McKean-Vlasov SDEs like Eq.(\ref{eq1}), say $\{X^{\e}\}$, $\e>0$ defined on $(\Omega, \mathscr{F}, \{\mathscr{F}_t\}_{t\in[0,T]}, \mP)$, which converge exponentially fast as $\e\rightarrow 0$. The decay rate is expressed via a rate function. An equivalent argument of the large deviation principle is the Laplace principle.

Next, we introduce the Laplace principle. Let us begin with some notations. Let $\mS$ be a Polish space. For each $\e>0$, let $X^{\e}$ be a $\mS$-valued random variable given on $(\Omega, \mathscr{F}, \{\mathscr{F}_t\}_{t\in[0,T]}, \mP)$.

\bd\label{compleve}
The function $I$ on $\mathbb{S}$ is called a rate function if for each $M<\infty$, $\{x\in \mathbb{S}:I(x)\leq M\}$ is a compact subset of $\mathbb{S}$.
\ed

\bd
We say that $\{X^{\epsilon}\}$ satisfies the Laplace principle with the rate function $I$,  if for and any real bounded continuous function $g$ on $\mathbb{S}$,
\ce
\lim\limits_{\varepsilon\rightarrow 0}\varepsilon \log \mE\left\{\exp\left[-\frac{g(X^{\epsilon})}{\epsilon}\right]\right\}=-\inf\limits_{x\in \mathbb{S}}\(g(x)+I(x)\).
\de
In particular, the family of $\{X^{\epsilon}, \e\in(0,1)\}$ satisfies the large deviation principle in $(\mS,\sB(\mS))$ with the rate function $I$. More precisely, for any closed subset $B_1\in \sB(\mS)$,
$$
\limsup\limits_{\e\rightarrow 0}\e\log\mP(X^{\epsilon}\in B_1)\leq -\inf\limits_{x\in B_1}I(x),
$$
and for any open subset $B_2\in \sB(\mS)$,
\ce
\liminf_{\e\rightarrow 0}\e\log\mP(X^{\epsilon}\in B_2)\geq -\inf\limits_{x\in B_2}I(x).
\de
\ed

To state some conditions under which the Laplace principle holds, we define some spaces. Set $\mathbb{H}:=L^{2}([0,T]; \mR^{m})$ and $\|h\|_{\mathbb{H}}:=(\int_{0}^{T}|h(t)|^{2}\dif t)^{\frac{1}{2}}$ for $h\in\mH$. Let $\mathcal{A}$ be the collection of predictable processes $u(\omega, \cdot)$ belonging to $\mathbb{H}$ a.s. $\omega$. For each $N\in\mN$ we define two following spaces:
\ce
\mathbf{D}_{2}^{N}:=\left\{h\in\mathbb{H}: \|h\|_{\mathbb{H}}^{2}\leq N \right\}, \quad \mathbf{A}_{2}^{N}:=\left\{u\in\mathcal{A}: u(\omega, \cdot)\in\mathbf{D}_{2}^{N}, a.s. \omega \right\}.
\de

\bco\label{cond}
Let $\psi^{\epsilon} : C([0,T];\mathbb{R}^m)\mapsto\mS$ be  a family of measurable mappings. There exists a
measurable mapping $\psi^{0} : C([0,T];\mathbb{R}^m)\mapsto\mS$ such that

$(i)$ for $N\in\mathbb{N}$ and $\{h_{\e}, \e>0\}\subset\mathbf{D}_2^{N}$, $h\in \mathbf{D}_2^{N}$, if $h_{\e}\rightarrow h$ as $\e\rightarrow 0$, then
\ce
\psi^{0}(\int_{0}^{\cdot}h_{\e}(s)\dif s)\longrightarrow \psi^{0}(\int_{0}^{\cdot}h(s)\dif s).
\de

$(ii)$ for $N\in\mathbb{N}$ and $\{u_{\epsilon},\epsilon>0\}\subset \mathbf{A}_{2}^{N}$, $u\in\mathbf{A}_{2}^{N}$, if $u_{\epsilon}$ converges in distribution to $u$ as $\e\rightarrow 0$, then 
$$
\psi^{\e}(\sqrt{\e}W(\cdot)+\int_{0}^{\cdot}u_{\e}(s)\dif s) \overset{d}{\longrightarrow} \psi^{0}(\int_{0}^{\cdot}u(s)\dif s).
$$
\eco

For $x\in\mS$ define ${\bf D}_{x}=\{h\in\mH: x=\psi^{0}(\int_{0}^{\cdot}h(s)\dif s)\}$. Let $I:\mathbb{S}\mapsto [0,\infty]$ be defined by               
$$
I(x)=\frac{1}{2}\inf\limits_{h\in{\bf D}_{x}}\|h\|_{\mathbb{H}}^{2}.
$$
From \cite[Theorem 4.4]{V}, we have the following result.

\bt\label{ldpbase}
Set $X^{\epsilon}:=\psi^{\epsilon}(\sqrt{\e}W)$. Assume that Condition \ref{cond} holds. Then $\{X^{\epsilon}\}$ satisfies the Laplace principle with the rate function $I$ given above.
\et

\subsection{The Laplace principle}

In the subsection, we study the Laplace principle for multivalued McKean-Vlasov SDEs.

Consider Eq.(\ref{meq1}), i.e.
\ce\left\{\begin{array}{l}
\dif X^{\epsilon}_{t}\in -A(X^{\epsilon}_{t})\dif t+b(X^{\epsilon}_{t},\sL_{X^{\epsilon}_{t}})\dif t+\sqrt{\epsilon}\sigma(X^{\epsilon}_{t},\sL_{X^{\epsilon}_{t}})\dif W_t, \\
 X^{\epsilon}_{0}=\xi\in\overline{\cD(A)}.
\end{array}
\right.
\de
Assume:
\begin{enumerate}[($\bf{H}_{1}$)]
	\item The function $b$ is continuous in $(x,\mu)$, and $b,\sigma$ satisfy for $(x,\mu)\in\mR^{d}\times{\cP_{2}(\mR^d)}$
	\ce
	|{b(x,\mu)}|\leq L_1(1+|x|+\|\mu\|_2), \quad \|\sigma(x,\mu)\|\leq L_1,
	\de
	where $L_1>0$ is a constant.
\end{enumerate}
\begin{enumerate}[($\bf{H}_{2}$)]
	\item The function $b,\sigma$ satisfy for $(x_{1},\mu_{1}),(x_{2},\mu_{2})\in\mR^{d}\times{\cP_{2}(\mR^d)}$
	\ce
	&&2\langle x_{1}-x_{2},b(x_{1},\mu_{1})-b(x_{2},\mu_{2})\rangle\leq L_2(|x_{1}-x_{2}|^{2}+\mW_2^{2}(\mu_{1},\mu_{2})),\\
	&&\|\sigma(x_{1},\mu_{1})-\sigma(x_{2},\mu_{2})\|^2\leq L_2(|x_{1}-x_{2}|^{2}+\mW_2^{2}(\mu_{1},\mu_{2})),
	\de
	where $L_2>0$ is a constant.
\end{enumerate}

Under $(\bf{H}_{1})$-$(\bf{H}_{2})$, by \cite[Theorem 3.5]{G}, we know that Eq.(\ref{meq1}) has a unique strong solution denoted as $(X^{\epsilon}, K^{\epsilon})$. In order to prove the Laplace principle for Eq.(\ref{eq1}),
we will verify Condition \ref{cond} with
\ce
\mathbb{S}:=C([0,T],\overline{\cD(A)}), \quad \psi^{\epsilon}(\sqrt{\e}W):=X^{\epsilon}.
\de

First of all, we consider a controlled  analogue of Eq.(\ref{meq1}) with the same initial value
\be
X^{\epsilon,u}_{t}&=&\xi+\int_{0}^{t}b(X^{\epsilon,u}_{s},\sL_{X^{\epsilon}_{s}})\dif s+\int_{0}^{t}\sigma(X^{\epsilon,u}_{s},\sL_{X^{\epsilon}_{s}})u(s)\dif s\no\\
&&+\sqrt{\epsilon}\int_{0}^{t}\sigma(X^{\epsilon,u}_{s},\sL_{X^{\epsilon}_{s}})\dif W_s- K^{\epsilon,u}_{t}, \quad u\in\mathbf{A}_2^{N}.
\label{contanal}
\ee
Thus, by the Girsanov theorem, it holds that Eq.(\ref{contanal}) has a unique solution denoted as $(X^{\epsilon,u}, K^{\epsilon,u})$. Moreover, $X^{\epsilon,u}=\psi^{\epsilon}(\sqrt{\epsilon}W+\int_{0}^{\cdot}u(s)\dif s)$.
Then let $ (X^{0,u}, K^{0,u})$ solve the following equation
\be
 X^{0,u}_{t}=\xi+\int_{0}^{t}b(X^{0,u}_{s}, \delta_{X^{0}_{s}})\dif s+\int_{0}^{t}\sigma(X^{0,u}_{s}, \delta_{X^{0}_{s}})u(s)\dif s- K^{0,u}_{t},
 \label{deteequa}
\ee
where $\delta_{\cdot}$ is the Dirac measure, i.e. for any $B\in\sB(\mR^d)$
$$
\delta_{x}(B)=\left\{\begin{aligned}
&1, x\in B,&\\
&0, x\notin B,&
\end{aligned}
\right.
$$
and $(X^{0}, K^0)$ is the unique solution of Eq.(\ref{eq3}), i.e.
\ce
X_t^{0}=\xi+\int_0^tb(X_s^{0},\delta_{X_s^{0}})\dif s-K_t^0.
\de
We define the measurable map $\psi^{0}: C([0,T]; \mR^{m})\mapsto \mS$ by $\psi^{0}(\int_{0}^{\cdot}u(s)\dif s)=X^{0,u}$. So, for $\psi^{0}$ we verify Condition \ref{cond}.

\bl\label{auxilemm3}
Let $h_{\e}\rightarrow h$ in $\mathbf{D}_2^{N}$ as ${\e}\rightarrow0$. Then $\psi^{0}(\int_{0}^{\cdot}h_{\e}(s)\dif s)$ converges to $\psi^{0}(\int_{0}^{\cdot}h(s)\dif s)$.
\el
\begin{proof}
First of all, by the definition of $\psi^{0}$, it holds that $\psi^{0}(\int_{0}^{\cdot}h_{\e}(s)\dif s), \psi^{0}(\int_{0}^{\cdot}h(s)\dif s)$ satisfy the following equations:
\ce
&&X^{0,h_{\e}}_{t}=\xi+\int_{0}^{t}b(X^{0,h_{\e}}_{s}, \delta_{X^{0}_{s}})\dif s+\int_{0}^{t}\sigma(X^{0,h_{\e}}_{s}, \delta_{X^{0}_{s}})h_{{\e}}(s)\dif s- {K}^{0,h_{\e}}_{t},\\
&&X^{0,h}_{t}=\xi+\int_{0}^{t}b(X^{0,h}_{s}, \delta_{X^{0}_{s}})\dif s+\int_{0}^{t}\sigma(X^{0,h}_{s}, \delta_{X^{0}_{s}})h(s)\dif s- {K}^{0,h}_{t}.
\de
Set $Z^{0}(t)=X^{0,h_{{\e}}}_{t}-X^{0,h}_{t}$, and by Lemma \ref{equi} and ($\bf{H}_{2}$) we have
\be
|Z^{0}(t)|^{2}&=&-2\int_{0}^{t}\langle Z^{0}(s),\dif {K}^{0,h_{{\e}}}_{s}-\dif {K}^{0,h}_{s}\rangle+2\int_{0}^{t}\langle Z^{0}(s),b(X^{0,h_{{\e}}}_{s},\delta_{X^{0}_{s}})-b(X^{0,h}_{s},\delta_{X^{0}_{s}})\rangle\dif s  \no\\
&&+2\int_{0}^{t}\langle Z^{0}(s),\sigma(X^{0,h_{{\e}}}_{s},\delta_{X^0_{s}})h_{{\e}}(s)-\sigma(X^{0,h}_{s},\delta_{X^0_{s}})h(s)\rangle\dif s  \no\\
&\leq&2\int_{0}^{t}\langle Z^{0}(s),b(X^{0,h_{{\e}}}_{s},\delta_{X^0_{s}})-b(X^{0,h}_{s},\delta_{X^0_{s}})\rangle\dif s  \no\\
&&+2\int_{0}^{t}\langle Z^{0}(s),\sigma(X^{0,h_{{\e}}}_{s},\delta_{X^0_{s}})h_{{\e}}(s)-\sigma(X^{0,h}_{s},\delta_{X^0_{s}})h(s)\rangle\dif s  \no\\
&\leq&L_{2}\int_{0}^{t}|Z^{0}(s)|^2\dif s+ 2\int_{0}^{t}\langle Z^{0}(s),\sigma(X^{0,h}_{s},\delta_{X^0_{s}})(h_{{\e}}(s)-h(s))\rangle \dif s\no\\
&&+ 2\int_{0}^{t}\langle Z^{0}(s),(\sigma(X^{0,h_{{\e}}}_{s},\delta_{X^0_{s}})- \sigma(X^{0,h}_{s},\delta_{X^0_{s}}))h_{{\e}}(s)\rangle \dif s\no\\
&=:&L_{2}\int_{0}^{t}|Z^{0}(s)|^2\dif s+I_1+I_2.
\label{eq10}
\ee

Next, for $I_1$, by H\"older's inequality and ($\bf{H}_{1}$), we have
\be
I_1&\leq&2\left|\int_{0}^{t}\langle Z^{0}(s),\sigma(X^{0,h}_{s}, \delta_{X^0_{s}})(h_{{\e}}(s)-h(s))\rangle \dif s\right| \no\\
&\leq& 2\sup\limits_{s\in[0,t]}|Z^{0}(s)|\int_{0}^{t}|\sigma(X^{0,h}_{s},\delta_{X^0_{s}})(h_{{\e}}(s)-h(s))|\dif s\no\\
&\leq& 2\sup\limits_{s\in[0,t]}|Z^{0}(s)|\left (\int_{0}^{t}\|\sigma(X^{0,h}_{s},\delta_{X^0_{s}})\|^{2}\dif s \right )^{\frac{1}{2}}\left(\int_{0}^{t}|h_{{\e}}(s)-h(s)|^{2}\dif s\right)^{\frac{1}{2}}\no\\
&\leq&\frac{1}{4}\sup\limits_{s\in[0,t]}|Z^{0}(s)|^{2}+C\int_{0}^{T}|h_{{\e}}(s)-h(s)|^{2}\dif s.
\label{fkf}
\ee
Besides, for $I_2$, from ($\bf{H}_{2}$), it follows that
\be
I_2&\leq&2\left|\int_{0}^{t}\langle Z^{0}(s),(\sigma(X^{0,h_{\e}}_{s},\delta_{X^0_{s}})- \sigma(X^{0,h}_{s}, \delta_{X^0_{s}}))h_{{\e}}(s)\rangle \dif s\right|\no\\
&\leq& 2\sqrt{L_{2}}\int_{0}^{t}|Z^{0}(s)|^{2}|h_{{\e}}(s)|\dif s\no\\
&\leq& 2\sqrt{L_{2}}\left(\int_{0}^{t}|Z^{0}(s)|^{4}\dif s\right)^{\frac{1}{2}}\left(\int_{0}^{t}|h_{{\e}}(s)|^{2}\dif s\right)^{\frac{1}{2}}\no\\
&\leq& 2\sqrt{L_{2}N}\left(\int_{0}^{t}|Z^{0}(s)|^{4}\dif s\right)^{\frac{1}{2}}\no\\
&\leq& 2\sqrt{L_{2}N}\left(\sup\limits_{s\in[0,t]}|Z^{0}(s)|\right)\left(\int_{0}^{t}|Z^{0}(s)|^{2}\dif s\right)^{\frac{1}{2}}\no\\
&\leq& \frac{1}{4}\sup\limits_{s\in[0,t]}|Z^{0}(s)|^{2}+C\int_{0}^{t}|Z^{0}(s)|^{2}\dif s.
\label{eq11}
\ee
Inserting (\ref{fkf}) and (\ref{eq11}) in (\ref{eq10}), we obtain that
\be
|Z^{0}(t)|^{2}\leq \frac{1}{2}\sup\limits_{s\in[0,t]}|Z^{0}(s)|^{2}+C\int_{0}^{t}|Z^{0}(s)|^{2}\dif s+C\int_{0}^{T}|h_{\e}(s)-h(s)|^{2}\dif s,
\ee
which implies that
$$
\sup\limits_{t\in[0,T]}|Z^{0}(t)|^{2}\leq C\int_{0}^{T}\left(\sup\limits_{s\in[0,u]}|Z^{0}(s)|^{2}\right)\dif u+C\int_{0}^{T}|h_{\e}(s)-h(s)|^{2}\dif s.
$$
Thus, by Gronwall's inequality, it holds that
\ce
\sup\limits_{t\in[0,T]}|X^{0,h_{{\e}}}_{t}-X^{0,h}_{t}|^2
\leq \left[C\int_{0}^{T}|h_{{\e}}(s)-h(s)|^{2}\dif s\right]e^{CT}.
\de
Therefore, we know that $\sup\limits_{t\in[0,T]}|X^{0,h_{{\e}}}_{t}-X^{0,h}_{t}|^2\rightarrow 0$ as ${\e}\rightarrow\infty$.
\end{proof}

\bl\label{xex0}
Under ($\bf{H}_{1}$)-($\bf{H}_{2}$), it holds that
$$
\mE\left(\sup\limits_{t\in[0,T]}|X_t^\e-X_t^0|^2\right)\leq C\e e^{CT}.
$$
\el
\begin{proof}
Note that
\ce
X_t^\e-X_t^0=-(K_t^\e-K_t^0)+\int_0^t(b(X_s^\e,\sL_{X_s^\e})-b(X_s^0,\d_{X_s^0}))\dif s+\sqrt{\epsilon}\int_{0}^{t}\sigma(X^{\epsilon}_{s},\sL_{X^{\epsilon}_{s}})\dif W_s.
\de
Thus, set $Z^\e(t):=X_t^\e-X_t^0$, and then it follows from the It\^o formula and Lemma \ref{equi}  that
\ce
|Z^\e(t)|^{2}&=&-2\int_{0}^{t}\langle Z^\e(s),\dif K^{\epsilon}_{s}- \dif K^{0,u}_{s}    \rangle+2\int_{0}^{t}\langle   Z^{\epsilon}(s), b(X^{\epsilon}_{s},\sL_{X^{\epsilon}_{s}})-b(X^{0}_{s},\delta_{X^0_{s}}) \rangle  \dif s   \no\\
&&+2\sqrt{\epsilon} \int_{0}^{t}\langle  Z^{\epsilon}(s),  \sigma (X^{\epsilon}_{s},\sL_{X^{\epsilon}_{s}})\dif W_s\rangle +\epsilon\int_{0}^{t}\|\sigma(X^{\epsilon}_{s},\sL_{X^{\epsilon}_{s}})\|^{2} \dif s\no\\
&\leq&2\int_{0}^{t}\langle   Z^{\epsilon}(s), b(X^{\epsilon}_{s},\sL_{X^{\epsilon}_{s}})-b(X^{0}_{s},\delta_{X^0_{s}}) \rangle  \dif s   \no\\
&&+2\sqrt{\epsilon} \int_{0}^{t}\langle  Z^{\epsilon}(s),  \sigma (X^{\epsilon}_{s},\sL_{X^{\epsilon}_{s}}) \dif W_s\rangle+\epsilon\int_{0}^{t}\|\sigma(X^{\epsilon}_{s},\sL_{X^{\epsilon}_{s}})\|^{2} \dif s\\
&=:&Q_1(t)+Q_2(t)+Q_3(t).
\de

For $Q_{1}(t)$, by ($\bf{H}_{2}$) we obtain that
\ce
Q_{1}(t)\leq L_2\int_{0}^{t}\left(|Z^{\epsilon}(s)|^2+\mW_2^2(\sL_{X^{\epsilon}_{s}},\delta_{X^0_{s}})\right)\dif s.
\de
Note that
\be
\mW_2^{2}(\sL_{X^{\epsilon}_{s}},\delta_{X^0_{s}})\leq\mE|X^{\epsilon}_{s}-X^0_{s}|^{2}.
\label{eq14}
\ee
So, one can get that
\be
\mE\left(\sup\limits_{t\in[0,T]}|Q_{1}(t)|\right)\leq 2L_2\mE\int_{0}^{T}| Z^{\epsilon}(s)|^{2} \dif s.
\label{q1}
\ee
For $Q_{2}(t)$, by the Burkholder-Davis-Gundy inequality and ($\bf{H}_{1}$), it holds that
\be
\mE\left(\sup\limits_{t\in[0,T]}|Q_{2}(t)|\right)&\leq& 2\sqrt{\epsilon}C\mE\left(\int_{0}^{T}|Z^{\epsilon}(s)|^2\|\sigma (X^{\epsilon}_{s},\sL_{X^{\epsilon}_{s}})\|^2 \dif s\right)^{1/2}\no\\
&\leq& 2\sqrt{\epsilon} L_1C\left(\int_{0}^{T}\mE|Z^{\epsilon}(s)|^2\dif s\right)^{1/2}\no\\
&\leq&C\epsilon+C\int_{0}^{T}\mE|Z^{\epsilon}(s)|^2\dif s.
\label{q2}
\ee
For $Q_{3}(t)$, from ($\bf{H}_{1}$), we have that
\be
\mE\left(\sup\limits_{t\in[0,T]}|Q_{3}(t)|\right)\leq L^2_1T\epsilon.
\label{q3}
\ee

Finally, by (\ref{q1})-(\ref{q3}), it holds that
\ce
\mE\left(\sup\limits_{t\in[0,T]}|Z^\e(t)|\right)\leq C\e+C\int_{0}^{T}\mE\left(\sup\limits_{u\in[0,s]}|Z^{\epsilon}(u)|^2\right)\dif s,
\de
which together with the Gronwall inequality yields the required estimate.
\end{proof}

\bl\label{auxilemm2}
Assume that for $\epsilon\in(0,1)$ and $\{u_{\epsilon}\}\subset\mathbf{A}_{2}^{N}$, $u\in\mathbf{A}_{2}^{N}$, $u_{\epsilon}$ converges to $u$ almost surely as $\epsilon\rightarrow 0$. Then $\psi^{\epsilon}(\sqrt{\epsilon}W+\int_{0}^{\cdot}u_{\epsilon}(s)\dif s)\rightarrow \psi^{0}(\int_{0}^{\cdot}u(s)\dif s)$ in probability.
\el
\begin{proof}
First of all, note that
$$
X^{\epsilon,u_{\epsilon}}=\psi^{\epsilon}(\sqrt{\epsilon}W+\int_{0}^{\cdot}u_{\epsilon}(s)\dif s), \quad X^{0,u}=\psi^{0}(\int_{0}^{\cdot}u(s)\dif s).
$$
To obtain $\psi^{\epsilon}(\sqrt{\epsilon}W+\int_{0}^{\cdot}u_{\epsilon}(s)\dif s){\rightarrow} \psi^{0}(\int_{0}^{\cdot}u(s)\dif s)$ in probability, we only need to prove $X^{\epsilon,u_{\epsilon}}-X^{0,u}{\rightarrow}0$ in probability.
Since the mean square convergence implies the convergence in probability, we estimate $\mE\left(\sup\limits_{t\in[0,T]}|X^{\epsilon,u_{\epsilon}}_{t}-X^{0,u}_{t}|^2\right)$.

Set $Z^{\epsilon,u_{\epsilon}}(t)=X^{\epsilon,u_{\epsilon}}_{t}-X^{0,u}_{t}$, and we have
\ce
Z^{\epsilon,u_{\epsilon}}(t)
&=&\int_{0}^{t}\left[b(X^{\epsilon,u_{\epsilon}}_{s},\sL_{X^{\epsilon}_{s}})-b(X^{0,u}_{s},\delta_{X^0_{s}})\right]\dif s\no\\
&&+\int_{0}^{t}\left[\sigma(X^{\epsilon,u_{\epsilon}}_{s},\sL_{X^{\epsilon}_{s}})u_{\epsilon}(s)-\sigma(X^{0,u}_{s}, \delta_{X^0_{s}})u(s)\right]\dif s\no\\
&&+\sqrt{\epsilon}\int_{0}^{t}\sigma(X^{\epsilon,u_{\epsilon}}_{s},\sL_{X^{\epsilon}_{s}})\dif W_s-(K^{\epsilon,u_{\epsilon}}_{t}-K^{0,u}_{t}).
\de
By the It\^o formula and Lemma \ref{equi}, it holds that
\be
|Z^{\epsilon,u_{\epsilon}}(t)|^{2}&=&-2\int_{0}^{t}\langle Z^{\epsilon,u_{\epsilon}}(s),\dif K^{\epsilon,u_{\epsilon}}_{s}- \dif K^{0,u}_{s}    \rangle  \no\\
&&+2\int_{0}^{t}\langle   Z^{\epsilon,u_{\epsilon}}(s), b(X^{\epsilon,u_{\epsilon}}_{s},\sL_{X^{\epsilon}_{s}})-b(X^{0,u}_{s},\delta_{X^0_{s}}) \rangle  \dif s   \no\\
&&+2\int_{0}^{t}\langle   Z^{\epsilon,u_{\epsilon}}(s), \sigma(X^{\epsilon,u_{\epsilon}}_{s},\sL_{X^{\epsilon}_{s}})u_{\epsilon}(s)-\sigma(X^{0,u}_{s},\delta_{X^0_{s}})u(s) \rangle  \dif s  \no\\
&&+2\sqrt{\epsilon} \int_{0}^{t}\langle  Z^{\epsilon,u_{\epsilon}}(s),  \sigma (X^{\epsilon,u_{\epsilon}}_{s},\sL_{X^{\epsilon}_{s}})\dif W_s\rangle  \no\\
&&+\epsilon\int_{0}^{t}\|\sigma(X^{\epsilon,u_{\epsilon}}_{s},\sL_{X^{\epsilon}_{s}})\|^{2} \dif s\no\\
&\leq&2\int_{0}^{t}\langle   Z^{\epsilon,u_{\epsilon}}(s), b(X^{\epsilon,u_{\epsilon}}_{s},\sL_{X^{\epsilon}_{s}})-b(X^{0,u}_{s},\delta_{X^0_{s}}) \rangle  \dif s   \no\\
&&+2\int_{0}^{t}\langle   Z^{\epsilon,u_{\epsilon}}(s), \sigma(X^{\epsilon,u_{\epsilon}}_{s},\sL_{X^{\epsilon}_{s}})u_{\epsilon}(s)-\sigma(X^{0,u}_{s},\delta_{X^0_{s}})u(s) \rangle  \dif s  \no\\
&&+2\sqrt{\epsilon} \int_{0}^{t}\langle  Z^{\epsilon,u_{\epsilon}}(s),  \sigma (X^{\epsilon,u_{\epsilon}}_{s},\sL_{X^{\epsilon}_{s}}) \dif W_s\rangle \no\\
&&+\epsilon\int_{0}^{t}\|\sigma(X^{\epsilon,u_{\epsilon}}_{s},\sL_{X^{\epsilon}_{s}})\|^{2} \dif s\no\\
&=:&J_{1}(t)+J_{2}(t)+J_{3}(t)+J_{4}(t).
\label{j1j2j3j4}
\ee

For $J_{1}(t)$, by ($\bf{H}_{2}$) one obtain that
\ce
J_{1}(t)\leq L_2\int_{0}^{t}\left(|Z^{\epsilon,u_{\epsilon}}(s)|^2+\mW_2^2(\sL_{X^{\epsilon}_{s}},\delta_{X^0_{s}})\right)\dif s.
\de
So, by (\ref{eq14}) we furthermore have
\be
\mE\left(\sup\limits_{t\in[0,T]}|J_{1}(t)|\right)\leq L_2\mE\int_{0}^{T}| Z^{\epsilon,u_{\epsilon}}(s)|^{2} \dif s+L_2\mE\int_{0}^{T}|X^{\epsilon}_{s}-X^0_{s}|^{2}\dif s.
\label{j1}
\ee

For $J_{2}$, we rewrite it as
\ce
J_{2}(t)&=&2\int_{0}^{t}\langle   Z^{\epsilon,u_{\epsilon}}(s), (\sigma(X^{\epsilon,u_{\epsilon}}_{s},\sL_{X^{\epsilon}_{s}})-\sigma(X^{0,u}_{s},\delta_{X^0_{s}}))u_{\epsilon}(s) \rangle  \dif s\no\\
&&+2\int_{0}^{t}\langle   Z^{\epsilon,u_{\epsilon}}(s), \sigma(X^{0,u}_{s},\delta_{X^0_{s}})(u_{\epsilon}(s)-u(s)) \rangle  \dif s.\no\\
&=&J_{21}(t)+J_{22}(t).
\de
For $J_{21}(t)$, by ($\bf{H}_{2}$), (\ref{eq14}) and the H\"older inequality, it holds that
\ce
&&\mE\left(\sup\limits_{t\in[0,T]}|J_{21}(t)|\right)\\
&\leq& 2\sqrt{L_{2}}\mE\int_{0}^{T}|Z^{\epsilon,u_{\epsilon}}(s)|\(|Z^{\epsilon,u_{\epsilon}}(s)|^2+\mE|X^{\epsilon}_{s}-X^0_{s}|^2\)^{1/2}|u_{\epsilon}(s)|\dif s\no\\
&\leq&\frac{1}{4}\mE\left[\sup\limits_{t\in[0,T]}|Z^{\epsilon,u_{\epsilon}}(s)|^{2}\right]+C\mE\left(\int_{0}^{T}\(|Z^{\epsilon,u_{\epsilon}}(s)|^2+\mE|X^{\epsilon}_{s}-X^0_{s}|^2\)^{1/2}|u_{\epsilon}(s)|\dif s\right)^2\\
&\leq& \frac{1}{4}\mE\left[\sup\limits_{t\in[0,T]}|Z^{\epsilon,u_{\epsilon}}(s)|^{2}\right]+C\mE\left(\int_{0}^{T}\(|Z^{\epsilon,u_{\epsilon}}(s)|^2+\mE|X^{\epsilon}_{s}-X^0_{s}|^{2}\)\dif s\right)\left(\int_{0}^{T}|u_{\epsilon}(s)|^{2}\dif s\right)\no\\
&\leq& \frac{1}{4}\mE\left[\sup\limits_{t\in[0,T]}|Z^{\epsilon,u_{\epsilon}}(s)|^{2}\right]+C\mE\left[\int_{0}^{T}|Z^{\epsilon,u_{\epsilon}}(s)|^{2}\dif s\right]+C\mE\int_{0}^{T}|X^{\epsilon}_{s}-X^0_{s}|^{2}\dif s.
\de
And for $J_{22}(t)$, by ($\bf{H}_{1}$) and the Young inequality, we get
\ce
&&\mE\left(\sup\limits_{t\in[0,T]}|J_{22}(t)|\right)\no\\
&\leq& 2\mE\left[\sup\limits_{t\in[0,T]}|Z^{\epsilon,u_{\epsilon}}(s)|\left(\int_{0}^{T}\|\sigma(X^{0,u}_{s},\delta_{X^0_{s}})\|^{2}\dif s \right)^{\frac{1}{2}}  \left( \int_{0}^{T} |u_{\epsilon}(s)-u(s)|^{2}\dif s \right)^{\frac{1}{2}}\right]\no\\
&\leq&\frac{1}{4}\mE\left[\sup\limits_{t\in[0,T]}|Z^{\epsilon,u_{\epsilon}}(s)|^{2}\right]+C\mE\left[\left(\int_{0}^{T}\|\sigma(X^{0,u}_{s},\delta_{X^0_{s}})\|^{2}\dif s \right)\left( \int_{0}^{T} |u_{\epsilon}(s)-u(s)|^{2}\dif s \right)\right]\no\\
&\leq&\frac{1}{4}\mE\left[\sup\limits_{t\in[0,T]}|Z^{\epsilon,u_{\epsilon}}(s)|^{2}\right]+C\int_{0}^{T} \mE|u_{\epsilon}(s)-u(s)|^{2}\dif s.
\de
Thus, one can obtain that
\be
\mE\left(\sup\limits_{t\in[0,T]}|J_{2}(t)|\right)&\leq&\frac{1}{2}\mE\left[\sup\limits_{t\in[0,T]}|Z^{\epsilon,u_{\epsilon}}(s)|^{2}\right]+C\mE\left[\int_{0}^{T}|Z^{\epsilon,u_{\epsilon}}(s)|^{2}\dif s\right]\no\\
&&+C\mE\int_{0}^{T}|X^{\epsilon}_{s}-X^0_{s}|^{2}\dif s+C\int_{0}^{T} \mE|u_{\epsilon}(s)-u(s)|^{2}\dif s.
\label{j2}
\ee

For $J_{3}(t)$, from the Burkholder-Davis-Gundy inequality and ($\bf{H}_{1}$), it follows that
\be
\mE\left(\sup\limits_{t\in[0,T]}|J_{3}(t)|\right)&\leq& 2\sqrt{\epsilon}C\mE\left(\int_{0}^{T}|Z^{\epsilon,u_{\epsilon}}(s)|^2\|\sigma (X^{\epsilon,u_{\epsilon}}_{s},\sL_{X^{\epsilon}_{s}})\|^2 \dif s\right)^{1/2}\no\\
&\leq& 2\sqrt{\epsilon L_1}C\left(\int_{0}^{T}\mE|Z^{\epsilon,u_{\epsilon}}(s)|^2\dif s\right)^{1/2}\no\\
&\leq&C\epsilon+C\int_{0}^{T}\mE|Z^{\epsilon,u_{\epsilon}}(s)|^2\dif s.
\label{j3}
\ee
For $J_{4}(t)$, by ($\bf{H}_{1}$), we know
\be
\mE\left(\sup\limits_{t\in[0,T]}|J_{4}(t)|\right)\leq L_1T\epsilon.
\label{j4}
\ee

Combining (\ref{j1})-(\ref{j4}) with (\ref{j1j2j3j4}), we can get
\ce
\mE\left(\sup\limits_{t\in[0,T]}|Z^{\epsilon,u_{\epsilon}}(t)|^{2}\right)&\leq& C\mE\left[\int_{0}^{T}|Z^{\epsilon,u_{\epsilon}}(s)|^{2}\dif s\right]+C\mE\int_{0}^{T}|X^{\epsilon}_{s}-X^0_{s}|^{2}\dif s\no\\
&&+C\int_{0}^{T} \mE|u_{\epsilon}(s)-u(s)|^{2}\dif s+(C+2L_1T)\epsilon,
\de
and furthermore
\ce
\mE\left(\sup\limits_{t\in[0,T]}|Z^{\epsilon,u_{\epsilon}}(t)|^{2}\right)&\leq& C\bigg[\mE\int_{0}^{T}|X^{\epsilon}_{s}-X^0_{s}|^{2}\dif s+\int_{0}^{T}\mE|u_{\epsilon}(t)-u(t)|^{2}\dif t+(C+2L_1T)\epsilon \bigg].
\de
As $\e\rightarrow 0$, the above inequality together with Lemma \ref{xex0} and the dominated convergence theorem implies that $X^{\epsilon,u_{\epsilon}}-X^{0,u}\rightarrow 0$ in the mean square.
The proof is complete.
\end{proof}

\bt\label{ldpmmv}
Assume that $(\bf{H}_{1})$ and $(\bf{H}_{2})$ hold. Then the family $\{X^{\epsilon},\epsilon\in[0,1]\}$ satisfies the large deviation principle in $\mS:=C([0,T],\overline{\mathcal{D}(A)})$ with the rate function given by
$$
I(x)=\frac{1}{2} \inf\limits_{h\in {\bf D}_{x}: x=X^{0,h}}\|h\|_{\mH}^2.
$$
\et
\begin{proof}
By Lemma \ref{auxilemm3}, we know that Condition \ref{cond} $(i)$ holds. 

Next, we verify Condition \ref{cond} $(ii)$. Then for $\epsilon\in(0,1)$ and $\{u_{\epsilon}, \e>0\}\subset\mathbf{A}_{2}^{N}$, $u\in\mathbf{A}_{2}^{N}$, let $u_{\epsilon}$ converge to $u$ in distribution. By the Skorohod theorem, there exists a 
probability space $(\tilde{\Omega}, \tilde{\sF}, \tilde{\mP})$, and ${\bf D}_2^N$-valued random variables $\{\tilde{u}_{\epsilon}\}$, $\tilde{u}$ and a $m$-dimensional Brownian motion $\tilde{W}$ defined on it such that

(i) $\sL_{(\tilde{u}_{\epsilon},\tilde{W})}=\sL_{(u_{\epsilon},W)}$ and $\sL_{\tilde{u}}=\sL_{u}$;

(ii) $\tilde{u}_{\epsilon}$ converges to $\tilde{u}$ almost surely.

In the following, we construct two multivalued McKean-Vlasov differential equations:
\be
X^{\epsilon,\tilde{u}_\e}_{t}&=&\xi+\int_{0}^{t}b(X^{\epsilon,\tilde{u}_\e}_{s},\sL_{X^{\epsilon}_{s}})\dif s+\int_{0}^{t}\sigma(X^{\epsilon,\tilde{u}_\e}_{s},\sL_{X^{\epsilon}_{s}})\tilde{u}_\e(s)\dif s\no\\
&&+\sqrt{\epsilon}\int_{0}^{t}\sigma(X^{\epsilon,\tilde{u}_\e}_{s},\sL_{X^{\epsilon}_{s}})\dif \tilde{W}_s-K^{\epsilon,\tilde{u}_\e}_{t},\label{contanal1}\\
 X^{0,\tilde{u}}_{t}&=&\xi+\int_{0}^{t}b(X^{0,\tilde{u}}_{s}, \delta_{X^0_{s}})\dif s+\int_{0}^{t}\sigma(X^{0,\tilde{u}}_{s}, \delta_{X^0_{s}})\tilde{u}(s)\dif s- K^{0,\tilde{u}}_{t}.
 \label{deteequa1}
\ee
Thus, Eq.(\ref{contanal1}) has a unique strong solution $X^{\epsilon,\tilde{u}_\e}$, and Eq.(\ref{deteequa1}) has a unique solution $X^{0,\tilde{u}}$. Moreover, it holds that
$$
X^{\epsilon,\tilde{u}_{\epsilon}}=\psi^{\epsilon}(\sqrt{\epsilon}\tilde{W}+\int_{0}^{\cdot}\tilde{u}_{\epsilon}(s)\dif s), \quad X^{0,\tilde{u}}=\psi^{0}(\int_{0}^{\cdot}\tilde{u}(s)\dif s).
$$
By Lemma \ref{auxilemm2}, we have that $\psi^{\epsilon}(\sqrt{\epsilon}\tilde{W}+\int_{0}^{\cdot}\tilde{u}_{\epsilon}(s)\dif s)\rightarrow \psi^{0}(\int_{0}^{\cdot}\tilde{u}(s)\dif s)$ in probability, which yields that
$\psi^{\epsilon}(\sqrt{\epsilon}\tilde{W}+\int_{0}^{\cdot}\tilde{u}_{\epsilon}(s)\dif s)\rightarrow \psi^{0}(\int_{0}^{\cdot}\tilde{u}(s)\dif s)$ in distribution. Note that 
\ce
\psi^{\epsilon}(\sqrt{\epsilon}\tilde{W}+\int_{0}^{\cdot}\tilde{u}_{\epsilon}(s)\dif s)&\overset{d}{=}&\psi^{\epsilon}(\sqrt{\epsilon}W+\int_{0}^{\cdot}u_{\epsilon}(s)\dif s),\\
\psi^{0}(\int_{0}^{\cdot}\tilde{u}(s)\dif s)&\overset{d}{=}&\psi^{0}(\int_{0}^{\cdot}u(s)\dif s).
\de
So, $\psi^{\epsilon}(\sqrt{\epsilon}W+\int_{0}^{\cdot}u_{\epsilon}(s)\dif s)\rightarrow \psi^{0}(\int_{0}^{\cdot}u(s)\dif s)$ in distribution. 

Finally, by Theorem \ref{ldpbase}, we draw the conclusion.
\end{proof}

\section{The central limit theorem for multivalued McKean-Vlasov SDEs}\label{clt}

In this section, we study the central limit theorem  for multivalued McKean-Vlasov SDEs. 

For $\epsilon>0$, consider the following multivalued McKean-Vlasov SDE:
\be\left\{\begin{array}{l}
\dif \frac{\hat{X}_t^{\epsilon}- X_t^{0}}{\sqrt{\epsilon}}\in \ -A(\frac{\hat{X}_t^{\epsilon}- X_t^{0}}{\sqrt{\epsilon}})\dif t+ \frac{b(\hat{X}_t^{\epsilon},\sL_{\hat{X}_t^{\epsilon}})-b(X_t^{0},\delta_{X_t^{0}})}{\sqrt{\epsilon}}\dif t+\sigma(\hat{X}_t^{\epsilon},\sL_{\hat{X}_t^{\epsilon}})\dif W_t, \quad t\in[0,T],\\
\frac{\hat{X}_0^{\epsilon}- X_0^{0}}{\sqrt{\epsilon}}=0, 
\end{array}
\label{eq2}
\right.
\ee
where $X^{0}$ satisfies Eq.(\ref{eq3}), i.e.
\ce\left\{\begin{array}{l}
\dif X_t^{0}\in \ -A(X_t^{0})\dif t+  b(X_t^{0},\delta_{X_t^{0}})\dif t,\\
X^0_0=\xi.
\end{array}
\right.
\de
Assume:
\begin{enumerate}[($\bf{H}_{3}$)]
	\item $b$ and $\sigma$ are continuous and  satisfy for $(x,\mu), (x_1,\mu_1), (x_2,\mu_2)\in\mR^{d}\times{\cP_{2}(\mR^d)}$:
	\ce
	&&\|\nabla b(x,\mu)\|\leq L_3, \quad \|D^{L}b(x,\mu)\|_{T_{\mu,2}} \leq L_3,\quad |b(0,\d_0)|\leq L_3,\\
    &&\|\sigma(x_1,\mu_1)-\sigma(x_2,\mu_2)\|\leq L_3(|x_1-x_2|+\mW_2(\mu_1,\mu_2)),\\
   &&\|\sigma(x,\mu)\|\leq L_3(1+|x|+\|\mu\|_2),
    \de
    where $L_3>0$ is a constant.
\end{enumerate}
\begin{enumerate}[($\bf{H}_{4}$)]
	\item For $(x_1,\mu_1),(x_2,\mu_2)\in\mR^{d}\times{\cP_{2}(\mR^d)}$ and $X$, $Y$, $\phi \in L^{2}(\Omega\mapsto\mR^d,\mP)$
	\ce
	&&\|\nabla b(x_1,\mu_1)-\nabla b(x_2,\mu_2)\|\leq L_4(|x_1-x_2|+\mW_2(\mu_1,\mu_2)),
	\de
    \ce
    &&|\mE\left\langle D^{L}b(x_1,\sL_{X})(X),\phi\right\rangle-\mE\left\langle D^{L}b(x_2,\sL_{Y})(Y),\phi\right\rangle|\no\\
    &\leq&L_4\left( |x_1-x_2|+\mW_2(\sL_{X},\sL_{Y})+(\mE|X-Y|^{2})^{\frac{1}{2}}   \right)\left( \mE|\phi|^{2}  \right)^{\frac{1}{2}},
    \de
    where $L_4>0$ is a constant.
\end{enumerate}

By ($\bf{H}_{3}$), it holds that for $(x_1,\mu_1), (x_2,\mu_2), (x,\mu)\in\mR^{d}\times{\cP_{2}(\mR^d)}$
\be
&&|b(x_1,\mu_1)-b(x_2,\mu_2)|\leq L_3(|x_1-x_2|+\mW_2(\mu_1,\mu_2)), \label{h31}\\
&&|b(x,\mu)|\leq L_3(1+|x|+\|\mu\|_2).
\label{h32}
\ee
Thus, under ($\bf{H}_{3}$), we know that Eq.(\ref{eq3}) has a unique solution $(X_{\cdot}^0, K_{\cdot}^0)$, and Eq.(\ref{eq2}) has a unique solution $(\frac{\hat{X}_{\cdot}^{\epsilon}- X_{\cdot}^{0}}{\sqrt{\epsilon}}, \hat{K}_{\cdot}^\e)$ (c.f. \cite[Theorem 3.5]{G}). Then we construct a multivalued McKean-Vlasov SDE:
\be\left\{\begin{array}{l}
\dif Z_{t}\in-A(Z_{t})\dif t+\nabla_{Z_{t}}b(X_{t}^{0},\delta_{X_{t}^{0}})\dif t+\mE\left\langle D^{L}b(X_{t}^{0},\delta_{X_{t}^{0}})(X_{t}^{0}),Z_{t}\right\rangle\dif t+\sigma(X_{t}^{0},\delta_{X_{t}^{0}})\dif W_{t},\\
Z_{0}=0,
\end{array}
\right.
\label{limiequa}
\ee
where $\nabla_{y}b(x,\mu)$ denotes the directional derivative of the function $b$ at $x$ in the direction $y$. And the assumption ($\bf{H}_{3}$) assures that Eq.(\ref{limiequa}) has a unique solution $(Z_{\cdot}, \hat{K}^0_{\cdot})$. So, the central limit theorem for Eq.(\ref{eq1}) means that
$$
\frac{\hat{X}_{\cdot}^{\epsilon}- X_{\cdot}^{0}}{\sqrt{\epsilon}}\overset{d}{\longrightarrow}Z_{\cdot}.
$$

Now we state the main result in the section.

\bt\label{cl}
Under assumptions $(\bf{H}_{3})$ and $(\bf{H}_{4})$, it holds that for $p\geq1$ 
\ce
\mE\left(\sup\limits_{t\in[0,T]}\left|\frac{\hat{X}_t^{\epsilon}- X_t^{0}}{\sqrt{\epsilon}}-Z_{t}\right|^{2p}\right)\leq C\epsilon^{p},
\de
where the constant $C>0$ is independent of $\e$.
\et

In order to prove the above theorem, we prepare some lemmas. The following result is from \cite{Re}.

\bl\label{auxilemm}
Let $(\Omega, \sF, \mathbb{P})$ be an atomless probability space, and let $X$, $Y\in L^{2}(\Omega\mapsto\mR^{d},\mP)$ with $\sL_{X}=\mu$. If either X and Y are bounded and f is L-differentiable at $\mu$, or $f\in C_{b}^{1}(\cP_2(\mR^d))$, then
$$
\lim\limits_{\epsilon\rightarrow0}\frac{f(\sL_{X+\epsilon Y})-f(\mu)}{\epsilon}=\mE \left\langle D^{L}f(\mu)(X),Y \right\rangle.
$$
Consequently,
$$
\left|\lim\limits_{\epsilon\rightarrow0}\frac{f(\sL_{X+\epsilon Y})-f(\mu)}{\epsilon}\right|=|\mE \left\langle D^{L}f(\mu)(X),Y \right\rangle|\leq\| D^{L}f(\mu) \|_{T_{\mu,2}}\sqrt{\mE|Y|^{2}}.
$$
\el

In the following we give some estimates.

\bl\label{auxilemm0}
Under the assumption $(\bf{H}_{3})$, it follows that for $p\geq1$
\ce
\sup\limits_{t\in[0,T]}|X_{t}^{0}|^{4p}\leq C.
\de
\el
\begin{proof}
Note that $X^0$ satisfies the following equation
$$
X_{t}^{0}=\xi-K^0_{t}+\int_{0}^{t}b(X_s^{0},\delta_{X_s^{0}})\dif s.
$$
Thus, for $\alpha\in Int(\cD(A))$, by the Taylor formula and Lemma \ref{inteineq}, it holds that
\ce
|X_{t}^{0}-\a|^{2}&=&|\xi-\a|^2+2\int_{0}^{t}\left\langle X_{s}^{0}-\a ,b(X_s^{0},\delta_{X_s^{0}})\right\rangle\dif s-2\int_{0}^{t}\left\langle X_{s}^{0}-\a,\dif K^0_{s}  \right\rangle\\
&\leq&|\xi-\a|^2+2\int_{0}^{t}\left\langle X_{s}^{0}-\a,b(X_s^{0},\delta_{X_s^{0}})\right\rangle\dif s+\gamma_{2}\int_{0}^{t} |X_{s}^0-\alpha|\dif s+\gamma_{3}t-\gamma_{1}|K^0|^t_{0}\\
&\leq&|\xi-\a|^2+(\g_2+\g_3)T+2\int_{0}^{t}\left\langle X_{s}^{0}-\a,b(X_s^{0},\delta_{X_s^{0}})\right\rangle\dif s+\gamma_{2}\int_{0}^{t} |X_{s}^0-\alpha|^2\dif s.
\de
By the H\"older inequality, the Young inequality and (\ref{h32}), we obtain that
\ce
\sup\limits_{s\in[0,t]}|X_{s}^{0}-\a|^{4p}&\leq&3^{2p-1}\(|\xi-\a|^2+(\g_2+\g_3)T\)^{2p}+3^{2p-1}\gamma^{2p}_{2}\left(\int_{0}^{t} |X_{u}^0-\alpha|^2\dif u\right)^{2p}\\
&&+3^{2p-1}2^{2p}\left(\int_{0}^{t}|\left\langle X_{u}^{0}-\a,b(X_u^{0},\delta_{X_u^{0}})\right\rangle|\dif u\right)^{2p}\\
&\leq&3^{2p-1}\(|\xi-\a|^2+(\g_2+\g_3)T\)^{2p}+3^{2p-1}\gamma^{2p}_{2}T^{2p-1}\int_{0}^{t} |X_{u}^0-\alpha|^{4p}\dif u\\
&&+3^{2p-1}2^{2p}T^{2p-1}\int_{0}^{t}|X_{u}^{0}-\a|^{2p}|b(X_u^{0},\delta_{X_u^{0}})|^{2p}\dif u\\
&\leq&3^{2p-1}\(|\xi-\a|^2+(\g_2+\g_3)T\)^{2p}+3^{2p-1}T^{2p-1}(2^{2p-1}+\gamma^{2p}_{2})\int_{0}^{t}|X_{u}^{0}-\a|^{4p}\dif u\\
&&+3^{2p-1}2^{2p-1}T^{2p-1}\int_{0}^{t}|b(X_u^{0},\delta_{X_u^{0}})|^{4p}\dif u\\
&\leq&3^{2p-1}\(|\xi-\a|^2+(\g_2+\g_3)T\)^{2p}+3^{2p-1}T^{2p-1}(2^{2p-1}+\gamma^{2p}_{2})\int_{0}^{t}|X_{u}^{0}-\a|^{4p}\dif u\\
&&+3^{2p-1}2^{2p-1}T^{2p-1}\int_{0}^{t}L_3^{4p}(1+2|X_u^{0}|)^{4p}\dif u\\
&\leq&C+C\int_{0}^{t}\sup\limits_{s\in[0,u]}|X_{s}^{0}-\a|^{4p}\dif u,
\de
which together with the Gronwall inequality yields the required result. The proof is complete.
\end{proof}

\bl\label{auxilemm1}
Under the assumption $(\bf{H}_{3})$, it holds that for $p\geq1$
\ce
\mE\left(\sup\limits_{t\in[0,T]}|Z_{t}^{\epsilon}|^{4p}\right)\leq C,\quad \mE\left(\sup\limits_{t\in[0,T]}|Z_{t}|^{4p}\right)\leq C,
\de
where $Z_{t}^{\epsilon}:=\frac{\hat{X}_t^{\epsilon}- X_t^{0}}{\sqrt{\epsilon}}$, and the constant $C>0$ is independent of $\epsilon$.
\el
\begin{proof}
Based on Eq.(\ref{eq2}), $Z_{\cdot}^{\epsilon}$ satisfies the following equation
$$
Z_{t}^{\epsilon}=\int_{0}^{t}\frac{b(\hat{X}_s^{\epsilon},\sL_{\hat{X}_s^{\epsilon}})-b(X_s^{0},\delta_{X_s^{0}})}{\sqrt{\epsilon}}\dif s+\int_{0}^{t}\sigma(\hat{X}_s^{\epsilon},\sL_{\hat{X}_s^{\epsilon}})\dif W_s-\hat{K}_{t}^{\epsilon}.
$$
For $\alpha\in Int(\cD(A))$, by the It\^o formula and Lemma \ref{inteineq}, we have that
\ce
|Z_{t}^{\epsilon}-\alpha|^{2} &=& |\alpha|^{2}+2\int_{0}^{t} \left\langle Z_{s}^{\epsilon}-\alpha,\frac{b(\hat{X}_s^{\epsilon},\sL_{\hat{X}_s^{\epsilon}})-b(X_s^{0},\delta_{X_s^{0}})}{\sqrt{\epsilon}} \right\rangle\dif s\no\\
&&+2\int_{0}^{t} \left\langle Z_{s}^{\epsilon}-\alpha, \sigma(\hat{X}_s^{\epsilon},\sL_{\hat{X}_s^{\epsilon}})\dif W_{s}\right\rangle\no\\
&&+\int_{0}^{t}\|\sigma(\hat{X}_s^{\epsilon},\sL_{\hat{X}_s^{\epsilon}} \|^{2}\dif s-2\int_{0}^{t} \left\langle Z_{s}^{\epsilon}-\alpha,\dif \hat{K}_{s}^{\epsilon} \right\rangle\no\\
&\leq&|\alpha|^{2}+2\int_{0}^{t} \left\langle Z_{s}^{\epsilon}-\alpha,\frac{b(\hat{X}_s^{\epsilon},\sL_{\hat{X}_s^{\epsilon}})-b(X_s^{0},\delta_{X_s^{0}})}{\sqrt{\epsilon}} \right\rangle\dif s\no\\
&&+2\int_{0}^{t} \left\langle Z_{s}^{\epsilon}-\alpha, \sigma(\hat{X}_s^{\epsilon},\sL_{\hat{X}_s^{\epsilon}})\dif W_{s}\right\rangle\no\\
&&+\int_{0}^{t}\|\sigma(\hat{X}_s^{\epsilon},\sL_{\hat{X}_s^{\epsilon}} \|^{2}\dif s+\gamma_{2}\int_{0}^{t} |Z_{s}^{\epsilon}-\alpha|\dif s+\gamma_{3}t-\gamma_{1}|\hat{K}^{\epsilon}|^t_{0}\no\\
&\leq&|\alpha|^{2}+2\int_{0}^{t} \left\langle Z_{s}^{\epsilon}-\alpha,\frac{b(\hat{X}_s^{\epsilon},\sL_{\hat{X}_s^{\epsilon}})-b(X_s^{0},\delta_{X_s^{0}})}{\sqrt{\epsilon}} \right\rangle\dif s\no\\
&&+2\int_{0}^{t} \left\langle Z_{s}^{\epsilon}-\alpha, \sigma(\hat{X}_s^{\epsilon},\sL_{\hat{X}_s^{\epsilon}})\dif W_{s}\right\rangle\no\\
&&+\int_{0}^{t}\|\sigma(\hat{X}_s^{\epsilon},\sL_{\hat{X}_s^{\epsilon}} \|^{2}\dif s+\gamma_{2}\int_{0}^{t} |Z_{s}^{\epsilon}-\alpha|^2\dif s+(\g_2+\gamma_{3})t.
\de
Moreover, it holds that
\be
\mE\left(\sup\limits_{s\in[0,t]}|Z_{s}^{\epsilon}-\alpha|^{4p}\right)&\leq& 5^{2p-1}\left(|\alpha|^{2}+(\g_2+\gamma_{3})T\right)^{2p}\no\\
&&+5^{2p-1}2^{2p}\mE \left[ \sup\limits_{s\in[0,t]}\left| \int_{0}^{s} \left\langle Z_{u}^{\epsilon}-\alpha,\frac{b(\hat{X}_u^{\epsilon},\sL_{\hat{X}_u^{\epsilon}})-b(X_u^{0},\delta_{X_u^{0}})}{\sqrt{\epsilon}} \right\rangle\dif u\right|^{2p} \right]\no\\
&&+5^{2p-1}2^{2p}\mE \left[ \sup\limits_{s\in[0,t]}\left| \int_{0}^{s} \left\langle Z_{u}^{\epsilon}-\alpha,\sigma(\hat{X}_u^{\epsilon},\sL_{\hat{X}_u^{\epsilon}})\dif W_{u} \right\rangle\right|^{2p} \right]\no\\
&&+5^{2p-1}\mE \left[ \sup\limits_{s\in[0,t]}\left(\int_{0}^{s}\|\sigma(\hat{X}_u^{\epsilon},\sL_{\hat{X}_u^{\epsilon}} \|^{2}\dif u\right)^{2p}  \right]\no\\
&&+5^{2p-1}\gamma^{2p}_{2}\mE \left[ \sup\limits_{s\in[0,t]}\left(\int_{0}^{s}|Z_{u}^{\epsilon}-\alpha|^2\dif u\right)^{2p} \right]\no\\
&=:&5^{2p-1}\left(|\alpha|^{2}+(\g_2+\gamma_{3})T\right)^{2p}+I_{1}+I_{2}+I_{3}+I_{4}.
\label{i1i2i3i4}
\ee

For $I_1$, from (\ref{h31}), the H\"older inequality and the Young inequality, it follows that
\be
I_{1}&\leq&5^{2p-1}2^{2p}\mE\left[\int_{0}^{t}|Z_{u}^{\epsilon}-\alpha|\left|\frac{b(\hat{X}_u^{\epsilon},\sL_{\hat{X}_u^{\epsilon}})-b(X_u^{0},\delta_{X_u^{0}})}{\sqrt{\epsilon}}\right| \dif u\right]^{2p}\no\\
&\leq&5^{2p-1}2^{2p}L_3^{2p}\mE\left[\int_{0}^{t}|Z_{u}^{\epsilon}-\alpha|\frac{|\hat{X}_u^{\epsilon}-X_u^{0}|+\mathbb{W}_{2}(\sL_{\hat{X}_u^{\epsilon}},\delta_{X_{u}^{0}})}{\sqrt{\epsilon}}\dif u\right]^{2p}\no\\
&\leq&5^{2p-1}2^{2p}L_3^{2p}\mE\left[\int_{0}^{t}|Z_{u}^{\epsilon}-\alpha|(|Z_{u}^{\epsilon}|+(\mE|Z_{u}^{\epsilon}|^{2})^{1/2})\dif u\right]^{2p}\no\\
&\leq&5^{2p-1}2^{2p}L_3^{2p}\mE\left[\int_{0}^{t}|Z_{u}^{\epsilon}-\alpha|^{2}\dif u\right]^{p}\left[\int_{0}^{t}(|Z_{u}^{\epsilon}|+(\mE|Z_{u}^{\epsilon}|^{2})^{1/2})^{2}\dif u\right]^{p}\no\\
&\leq&C\mE\left[\int_{0}^{t}|Z_{u}^{\epsilon}-\alpha|^{2}\dif u\right]^{2p}+C\mE\left[\int_{0}^{t}(|Z_{u}^{\epsilon}|^{2}+\mE|Z_{u}^{\epsilon}|^{2})\dif u\right]^{2p}\no\\
&\leq&C\mE\int_{0}^{t}|Z_{u}^{\epsilon}-\alpha|^{4p}\dif u+C\mE\int_{0}^{t}(|Z_{u}^{\epsilon}|^{4p}+\mE|Z_{u}^{\epsilon}|^{4p})\dif u\no\\
&\leq&C\int_{0}^{t}\mE\left(\sup\limits_{s\in[0,u]}|Z_{s}^{\epsilon}-\alpha|^{4p}\right)\dif u+C|\a|^{4p},
\label{i1}
\ee
where we use the fact that $\mathbb{W}_{2}(\sL_{\hat{X}_u^{\epsilon}},\delta_{X_{u}^{0}})\leq (\mE|\hat{X}_u^{\epsilon}-X_{u}^{0}|^{2})^{1/2}$.

For $I_2$, by the BDG inequality, the H\"older inequality, the Young inequality, Lemma \ref{auxilemm0} and ($\bf{H}_{3}$), we have
\be
I_{2}&\leq&5^{2p-1}2^{2p}C\mE \left[\int_{0}^{t}|Z_{u}^{\epsilon}-\alpha|^{2}\|\sigma(\hat{X}_u^{\epsilon},\sL_{\hat{X}_u^{\epsilon}}) \|^{2}\dif u\right]^{p}\no\\
&\leq&5^{2p-1}2^{2p}T^{p-1}C\mE\int_{0}^{t}|Z_{u}^{\epsilon}-\alpha|^{2p}\|\sigma(\hat{X}_u^{\epsilon},\sL_{\hat{X}_u^{\epsilon}}) \|^{2p}\dif u\no\\
&\leq&C\mE\int_{0}^{t}|Z_{u}^{\epsilon}-\alpha|^{4p}\dif u+C\mE\int_{0}^{t}\|\sigma(\hat{X}_u^{\epsilon},\sL_{\hat{X}_u^{\epsilon}}) \|^{4p}\dif u\no\\
&\leq&C\mE\int_{0}^{t}|Z_{u}^{\epsilon}-\alpha|^{4p}\dif u+C\mE\int_{0}^{t}(1+|\hat{X}_u^{\epsilon}|+\|\sL_{\hat{X}_u^{\epsilon}}\|_2)^{4p}\dif u\no\\
&\leq&C\mE\int_{0}^{t}|Z_{u}^{\epsilon}-\alpha|^{4p}\dif u+C+C\mE\int_{0}^{t}|\hat{X}_u^{\epsilon}|^{4p}\dif u\no\\
&\leq&C\mE\int_{0}^{t}|Z_{u}^{\epsilon}-\alpha|^{4p}\dif u+C+C\mE\int_{0}^{t}(|\sqrt{\e}Z_{u}^{\epsilon}|^{4p}+|X^0_u|^{4p})\dif u\no\\
&\leq&C\int_{0}^{t}\mE\left(\sup\limits_{s\in[0,u]}|Z_{s}^{\epsilon}-\alpha|^{4p}\right)\dif u+C.
\label{i2}
\ee

By the H\"older inequality, Lemma \ref{auxilemm0} and ($\bf{H}_{3}$), one can obtain that
\be
I_{3}+I_4 \leq C+C\int_{0}^{t}\mE\left(\sup\limits_{s\in[0,u]}|Z_{s}^{\epsilon}-\alpha|^{4p}\right)\dif u.
\label{i3i4}
\ee

Therefore, by (\ref{i1i2i3i4})-(\ref{i3i4}), we see that
\ce
\mE\left(\sup\limits_{s\in[0,t]}|Z_{s}^{\epsilon}-\alpha|^{4p}\right)\leq C+C\int_{0}^{t}\mE\left(\sup\limits_{s\in[0,u]}|Z_{s}^{\epsilon}-\alpha|^{4p}\right)\dif u,
\de
which together with the Gronwall inequality yields the first estimate.

Next, we compute the second estimate. Note that $Z_{t}$ satisfies the following equation
\ce
Z_{t}=-\hat{K}^0_{t}+\int_{0}^{t}\nabla_{Z_{s}}b(X_{s}^{0},\delta_{X_{s}^{0}})\dif s+\int_{0}^{t}\mE\left\langle D^{L}b(X_{s}^{0},\delta_{X_{s}^{0}})(X_{s}^{0}),Z_{s}\right\rangle\dif s+\int_{0}^{t}\sigma(X_{s}^{0},\delta_{X_{s}^{0}})\dif W_{s}.
\de
Thus, for $\alpha\in Int(\cD(A))$, by the It\^o formula, Lemma \ref{inteineq} and $(\bf{H}_{3})$, it holds that
\ce
|Z_{t}-\a|^{2}&=&|\a|^2-2\int_{0}^{t}\left\langle  Z_{s}-\a, \dif\hat{K}^0_{s}  \right\rangle+2\int_{0}^{t}\left\langle Z_{s}-\a, \nabla_{Z_{s}}b(X_{s}^{0},\delta_{X_{s}^{0}})  \right\rangle \dif s\no\\
&&+2\int_{0}^{t}\bigg\langle Z_{s}-\a,\mE\left\langle  D^{L}b(X_{s}^{0},\delta_{X_{s}^{0}})(X_{s}^{0}),Z_{s}\right\rangle   \bigg\rangle \dif s\no\\
&&+2\int_{0}^{t}\left\langle  Z_{s}-\a, \sigma(X_{s}^{0},\delta_{X_{s}^{0}})  \dif W_{s}\right\rangle+\int_{0}^{t} \|\sigma(X_{s}^{0},\delta_{X_{s}^{0}})\|^{2}   \dif s\no\\
&\leq&|\a|^2+\gamma_{2}\int_{0}^{t} |Z_{s}-\alpha|\dif s+\gamma_{3}t-\gamma_{1}|\hat{K}^0|^t_{0}+2\int_{0}^{t}\left\langle Z_{s}-\a, \nabla_{Z_{s}}b(X_{s}^{0},\delta_{X_{s}^{0}})  \right\rangle \dif s\\
&&+2\int_{0}^{t}\bigg\langle Z_{s}-\a,\mE\left\langle  D^{L}b(X_{s}^{0},\delta_{X_{s}^{0}})(X_{s}^{0}),Z_{s}\right\rangle   \bigg\rangle \dif s\no\\
&&+2\int_{0}^{t}\left\langle  Z_{s}-\a, \sigma(X_{s}^{0},\delta_{X_{s}^{0}})  \dif W_{s}\right\rangle+\int_{0}^{t} \|\sigma(X_{s}^{0},\delta_{X_{s}^{0}})\|^{2}   \dif s\\
&\leq&|\a|^2+(\g_2+\gamma_{3}+L_3^2)t+\gamma_{2}\int_{0}^{t} |Z_{s}-\alpha|^2\dif s+2\int_{0}^{t}\left\langle Z_{s}-\a, \nabla_{Z_{s}}b(X_{s}^{0},\delta_{X_{s}^{0}})  \right\rangle \dif s\\
&&+2\int_{0}^{t}\bigg\langle Z_{s}-\a,\mE\left\langle  D^{L}b(X_{s}^{0},\delta_{X_{s}^{0}})(X_{s}^{0}),Z_{s}\right\rangle   \bigg\rangle \dif s\no\\
&&+2\int_{0}^{t}\left\langle  Z_{s}-\a, \sigma(X_{s}^{0},\delta_{X_{s}^{0}})  \dif W_{s}\right\rangle.
\de
Moreover, we know that
\be
&&\mE\left(\sup\limits_{s\in[0,t]}|Z_{t}-\a|^{4p}\right)\no\\
&\leq&5^{2p-1}\left(|\a|^2+(\g_2+\gamma_{3}+L_3^2)T\right)^{2p}+5^{2p-1}\gamma^{2p}_{2}T^{2p-1}\int_{0}^{t} \mE\left(\sup\limits_{s\in[0,u]}|Z_{s}-\a|^{4p}\right)\dif u\no\\
&&+5^{2p-1}2^{2p}\mE\left[ \sup\limits_{s\in[0,t]}\bigg|\int_{0}^{s}\left\langle Z_{u}-\a,  \nabla_{Z_{u}}b(X_{u}^{0},\delta_{X_{u}^{0}})     \right\rangle\dif u\bigg|^{2p} \right]\no\\
&&+5^{2p-1}2^{2p}\mE\left[ \sup\limits_{s\in[0,t]}\left|\int_{0}^{s}\bigg\langle Z_{u}-\a, \mE\left\langle  D^{L}b(X_{u}^{0},\delta_{X_{u}^{0}})(X_{u}^{0}),Z_{u}\right\rangle     \bigg\rangle\dif u\right|^{2p} \right]\no\\
&&+5^{2p-1}2^{2p}\mE\left[ \sup\limits_{s\in[0,t]}\left|\int_{0}^{s}\bigg\langle Z_{u}-\a,  \sigma(X_{u}^{0},\delta_{X_{u}^{0}})       \dif W_{u}    \bigg\rangle\right|^{2p} \right]\no\\
&=:&5^{2p-1}\left(|\a|^2+(\g_2+\gamma_{3}+L_3^2)T\right)^{2p}+5^{2p-1}\gamma^{2p}_{2}T^{2p-1}\int_{0}^{t} \mE\left(\sup\limits_{s\in[0,u]}|Z_{s}-\a|^{4p}\right)\dif u\no\\
&&+J_{1}+J_{2}+J_{3}.
\label{j1j2j3}
\ee

By the Young inequality, the H\"older inequality and $(\bf{H}_{3})$, we deduce 
\be
J_{1}&\leq&5^{2p-1}2^{2p}T^{2p-1}\mE\int_{0}^{t}|Z_{u}-\a|^{2p}\left|\nabla_{Z_{u}}b(X_{u}^{0},\delta_{X_{u}^{0}}) \right|^{2p}\dif u\no\\
&\leq&5^{2p-1}2^{2p}T^{2p-1}\mE\int_{0}^{t}|Z_{u}-\a|^{2p}\|\nabla b(X_{u}^{0},\delta_{X_{u}^{0}})\|^{2p}|Z_{u}|^{2p}\dif u\no\\
&\leq&C\mE\int_{0}^{t}|Z_{u}-\a|^{2p}(|Z_{u}-\a|^{2p}+|\a|^{2p})\dif u\no\\
&=&C\mE\int_{0}^{t}|Z_{u}-\a|^{4p}\dif u+C|\a|^{2p}\mE\int_{0}^{t}|Z_{u}-\a|^{2p}\dif u\no\\
&\leq&C\int_{0}^{t}\mE\left(\sup\limits_{s\in[0,u]}|Z_{s}-\a|^{4p}\right)\dif u+C,
\label{j12}
\ee
and 
\be
J_{2}&\leq&5^{2p-1}2^{2p}T^{2p-1}\mE\int_{0}^{t}|Z_{u}-\a|^{2p}|\mE\left\langle  D^{L}b(X_{u}^{0},\delta_{X_{u}^{0}})(X_{u}^{0}),Z_{u}\right\rangle |^{2p}\dif u\no\\
&\leq&5^{2p-1}2^{2p}T^{2p-1}\mE\int_{0}^{t}|Z_{u}-\a|^{2p}\left(\mE\|D^{L}b(X_{u}^{0},\delta_{X_{u}^{0}})\|_{T_{\delta_{X_{u}^{0}},2}}|Z_{u}|\right)^{2p}\dif u\no\\
&\leq&C\mE\int_{0}^{t}|Z_{u}-\a|^{2p}\mE(|Z_{u}-\a|^{2p}+|\a|^{2p})\dif u\no\\
&\leq&C\int_{0}^{t}\mE\left(\sup\limits_{s\in[0,u]}|Z_{s}-\a|^{4p}\right)\dif u+C.
\label{j22}
\ee

From the BDG inequality, the Young inequality, the H\"older inequality and $(\bf{H}_{3})$, it follows that
\be
J_{3}&\leq&C\mE\bigg[\int_{0}^{t}|Z_{u}-\a|^{2}   \| \sigma(X_{u}^{0},\delta_{X_{u}^{0}}) \|^{2}  \dif u\bigg]^{p}\no\\
&\leq&C\mE\int_{0}^{t}|Z_{u}-\a|^{2p} \dif u\leq C\mE\int_{0}^{t}|Z_{u}-\a|^{4p} \dif u+C\no\\
&\leq&C\int_{0}^{t}\mE\left(\sup\limits_{s\in[0,u]}|Z_{s}-\a|^{4p}\right)\dif u+C.
\label{j32}
\ee

Finally, combining (\ref{j12})-(\ref{j32}) with (\ref{j1j2j3}), we obtain that
\ce
\mE\left(\sup\limits_{s\in[0,t]}|Z_{t}-\a|^{4p}\right)\leq C\int_{0}^{t}\mE\left(\sup\limits_{s\in[0,u]}|Z_{s}-\a|^{4p}\right)\dif u+C,
\de
which together with the Gronwall inequality implies the second estimate. The proof is complete.
\end{proof}

{\bf The proof of Theorem \ref{cl}}.

By the definitions of $Z_{t}^{\epsilon}$ and $Z_{t}$, we know that
\ce
&&Z_{t}^{\epsilon}-Z_{t}\no\\
&=&\int_{0}^{t}\left[\frac{b(\hat{X}_s^{\epsilon},\sL_{\hat{X}_s^{\epsilon}})-b(X_s^{0},\delta_{X_s^{0}})}{\sqrt{\epsilon}}-\left(\nabla_{Z_{s}}b(X_{s}^{0},\delta_{X_{s}^{0}})+\mE\left\langle D^{L}b(X_{s}^{0},\delta_{X_{s}^{0}})(X_{s}^{0}),Z_{s}\right\rangle \right)\right]\dif s\no\\
&&+\int_{0}^{t}\left(\sigma(\hat{X}_s^{\epsilon},\sL_{\hat{X}_s^{\epsilon}})-\sigma(X_{s}^{0},\delta_{X_{s}^{0}})\right)\dif W_s-(\hat{K}_{t}^{\epsilon}-\hat{K}^0_{t}),
\de
which together with the It\^o formula yields that
\ce
&&|Z_{t}^{\epsilon}-Z_{t}|^{2p}\no\\
&=&2p\int_{0}^{t}|Z_{s}^{\epsilon}-Z_{s}|^{2p-2}\left\langle Z_{s}^{\epsilon}-Z_{s},  \frac{b(\hat{X}_s^{\epsilon},\sL_{\hat{X}_s^{\epsilon}})-b(X_s^{0},\delta_{X_s^{0}})}{\sqrt{\epsilon}}\right\rangle\dif s\no\\
&&-2p\int_{0}^{t}|Z_{s}^{\epsilon}-Z_{s}|^{2p-2}\left\langle Z_{s}^{\epsilon}-Z_{s},   \nabla_{Z_{s}}b(X_{s}^{0},\d_{X_s^0})+\mE\left\langle D^{L}b(X_{s}^{0},\delta_{X_{s}^{0}})(X_{s}^{0}),Z_{s}\right\rangle       \right\rangle\dif s\no\\
&&+2p\int_{0}^{t}|Z_{s}^{\epsilon}-Z_{s}|^{2p-2}\left\langle Z_{s}^{\epsilon}-Z_{s},\left(   \sigma(\hat{X}_s^{\epsilon},\sL_{\hat{X}_s^{\epsilon}})-\sigma(X_{s}^{0},\delta_{X_{s}^{0}}) \right)\dif W_{s} \right\rangle\no\\
&&+2p(p-1)\int_{0}^{t}|Z_{s}^{\epsilon}-Z_{s}|^{2p-4}\Big\langle  Z_{s}^{\epsilon}-Z_{s},\left( \sigma(\hat{X}_s^{\epsilon},\sL_{\hat{X}_s^{\epsilon}})-\sigma(X_{s}^{0},\delta_{X_{s}^{0}})\right)\\
&&\qquad\qquad\times \left(\sigma(\hat{X}_s^{\epsilon},\sL_{\hat{X}_s^{\epsilon}})-\sigma(X_{s}^{0},\delta_{X_{s}^{0}}) \right)^{*}(Z_{s}^{\epsilon}-Z_{s})\Big\rangle\dif s\no\\
&&+p\int_{0}^{t}|Z_{s}^{\epsilon}-Z_{s}|^{2p-2}\|\sigma(\hat{X}_s^{\epsilon},\sL_{\hat{X}_s^{\epsilon}})-\sigma(X_{s}^{0},\delta_{X_{s}^{0}})\|^{2}\dif s\\
&&-2p\int_{0}^{t}|Z_{s}^{\epsilon}-Z_{s}|^{2p-2}\left\langle Z_{s}^{\epsilon}-Z_{s},\dif \hat{K}_{s}^{\epsilon}-\dif \hat{K}^0_{s} \right\rangle.
\de
By Lemma \ref{equi}, one can obtain that
\ce
&&|Z_{t}^{\epsilon}-Z_{t}|^{2p}\no\\
&\leq&2p\int_{0}^{t}|Z_{s}^{\epsilon}-Z_{s}|^{2p-2}\left\langle Z_{s}^{\epsilon}-Z_{s},  \frac{b(\hat{X}_s^{\epsilon},\sL_{\hat{X}_s^{\epsilon}})-b(X_s^{0},\delta_{X_s^{0}})}{\sqrt{\epsilon}}\right\rangle\dif s\no\\
&&-2p\int_{0}^{t}|Z_{s}^{\epsilon}-Z_{s}|^{2p-2}\left\langle Z_{s}^{\epsilon}-Z_{s},   \nabla_{Z_{s}}b(X_{s}^{0},\d_{X_s^0})+\mE\left\langle D^{L}b(X_{s}^{0},\delta_{X_{s}^{0}})(X_{s}^{0}),Z_{s}\right\rangle       \right\rangle\dif s\no\\
&&+2p\int_{0}^{t}|Z_{s}^{\epsilon}-Z_{s}|^{2p-2}\left\langle Z_{s}^{\epsilon}-Z_{s},\left(   \sigma(\hat{X}_s^{\epsilon},\sL_{\hat{X}_s^{\epsilon}})-\sigma(X_{s}^{0},\delta_{X_{s}^{0}}) \right)\dif W_{s} \right\rangle\no\\
&&+p(2p-1)\int_{0}^{t}|Z_{s}^{\epsilon}-Z_{s}|^{2p-2}\|\sigma(\hat{X}_s^{\epsilon},\sL_{\hat{X}_s^{\epsilon}})-\sigma(X_{s}^{0},\delta_{X_{s}^{0}})\|^{2}\dif s\\
&=&2p\int_{0}^{t}|Z_{s}^{\epsilon}-Z_{s}|^{2p-2}\left\langle Z_{s}^{\epsilon}-Z_{s}, \frac{b(\hat{X}_s^{\epsilon},\sL_{\hat{X}_s^{\epsilon}})-b(X_s^{0},\sL_{\hat{X}_s^{\epsilon}})}{\sqrt{\epsilon}}-  \nabla_{Z^\e_{s}}b(X_{s}^{0},\sL_{\hat{X}_s^{\epsilon}})\right\rangle\dif s\\
&&+2p\int_{0}^{t}|Z_{s}^{\epsilon}-Z_{s}|^{2p-2}\left\langle Z_{s}^{\epsilon}-Z_{s},\frac{b(X_s^{0},\sL_{\hat{X}_s^{\epsilon}})-b(X_s^{0},\delta_{X_s^{0}})}{\sqrt{\epsilon}}-   \mE\left\langle D^{L}b(X_{s}^{0},\delta_{X_{s}^{0}})(X_{s}^{0}),Z_{s}^{\epsilon}\right\rangle \right\rangle\dif s\\
&&+2p\int_{0}^{t}|Z_{s}^{\epsilon}-Z_{s}|^{2p-2}\left\langle Z_{s}^{\epsilon}-Z_{s},\nabla_{Z_{s}^{\epsilon}}b(X_{s}^{0},\sL_{\hat{X}_{s}^{\epsilon}})-\nabla_{Z_{s}}b(X_{s}^{0},\sL_{\hat{X}_{s}^{\epsilon}})\right\rangle\dif s\\
&&+2p\int_{0}^{t}|Z_{s}^{\epsilon}-Z_{s}|^{2p-2}\left\langle Z_{s}^{\epsilon}-Z_{s},\nabla_{Z_{s}}b(X_{s}^{0},\sL_{\hat{X}_{s}^{\epsilon}})-\nabla_{Z_{s}}b(X_{s}^{0},\delta_{X_{s}^{0}}) \right\rangle\dif s\\
&&+2p\int_{0}^{t}|Z_{s}^{\epsilon}-Z_{s}|^{2p-2}\left\langle Z_{s}^{\epsilon}-Z_{s},\mE\left\langle D^{L}b(X_{s}^{0},\delta_{X_{s}^{0}})(X_{s}^{0}),Z_{s}^{\epsilon}\right\rangle
-\mE\left\langle D^{L}b(X_{s}^{0},\delta_{X_{s}^{0}})(X_{s}^{0}),Z_{s}\right\rangle\right\rangle\dif s\\
&&+2p\int_{0}^{t}|Z_{s}^{\epsilon}-Z_{s}|^{2p-2}\left\langle Z_{s}^{\epsilon}-Z_{s},\left(   \sigma(\hat{X}_s^{\epsilon},\sL_{\hat{X}_s^{\epsilon}})-\sigma(X_{s}^{0},\delta_{X_{s}^{0}}) \right)\dif W_{s} \right\rangle\no\\
&&+p(2p-1)\int_{0}^{t}|Z_{s}^{\epsilon}-Z_{s}|^{2p-2}\|\sigma(\hat{X}_s^{\epsilon},\sL_{\hat{X}_s^{\epsilon}})-\sigma(X_{s}^{0},\delta_{X_{s}^{0}})\|^{2}\dif s.
\de
From the above inequality, it follows that
\be
\mE\left(\sup\limits_{s\in[0,t]}|Z_{s}^{\epsilon}-Z_{s}|^{2p}\right)\leq U_{1}+U_{2}+U_{3}+U_{4}+U_{5}+U_{6}+U_{7},
\label{allesti}
\ee
where 
\ce
&&U_1:=2p\mE\left(\sup\limits_{s\in[0,t]}\int_{0}^{s}|Z_{u}^{\epsilon}-Z_{u}|^{2p-1}\left|\frac{b(\hat{X}_u^{\epsilon},\sL_{\hat{X}_u^{\epsilon}})-b(X_u^{0},\sL_{\hat{X}_u^{\epsilon}})}{\sqrt{\epsilon}}-  \nabla_{Z^\e_{u}}b(X_{u}^{0},\sL_{\hat{X}_u^{\epsilon}})\right|\dif u\right),\\
&&U_2:=2p\mE\bigg(\sup\limits_{s\in[0,t]}\int_{0}^{s}|Z_{u}^{\epsilon}-Z_{u}|^{2p-1}\left|\frac{b(X_u^{0},\sL_{\hat{X}_u^{\epsilon}})-b(X_u^{0},\delta_{X_u^{0}})}{\sqrt{\epsilon}}-\mE\left\langle D^{L}b(X_{u}^{0},\delta_{X_{u}^{0}})(X_{u}^{0}),Z_{u}^{\epsilon}\right\rangle\right|\dif u\bigg),\\
&&U_3:=2p\mE\left(\sup\limits_{s\in[0,t]}\int_{0}^{s}|Z_{u}^{\epsilon}-Z_{u}|^{2p-1}\left|\nabla_{Z_{u}^{\epsilon}}b(X_{u}^{0},\sL_{\hat{X}_{u}^{\epsilon}})-\nabla_{Z_{u}}b(X_{u}^{0},\sL_{\hat{X}_{u}^{\epsilon}})\right|\dif u\right),\\
&&U_4:=2p\mE\left(\sup\limits_{s\in[0,t]}\int_{0}^{s}|Z_{u}^{\epsilon}-Z_{u}|^{2p-1}\left|\nabla_{Z_{u}}b(X_{u}^{0},\sL_{\hat{X}_{u}^{\epsilon}})-\nabla_{Z_{u}}b(X_{u}^{0},\delta_{X_{u}^{0}}) \right|\dif u\right),\\
&&U_5:=2p\mE\bigg(\sup\limits_{s\in[0,t]}\int_{0}^{s}|Z_{u}^{\epsilon}-Z_{u}|^{2p-1}|\mE\left\langle D^{L}b(X_{u}^{0},\delta_{X_{u}^{0}})(X_{u}^{0}),Z_{u}^{\epsilon}\right\rangle-\mE\left\langle D^{L}b(X_{u}^{0},\delta_{X_{u}^{0}})(X_{u}^{0}),Z_{u}\right\rangle|\dif u\bigg),\\
&&U_6:=2p\mE\left(\sup\limits_{s\in[0,t]}\int_{0}^{s}|Z_{u}^{\epsilon}-Z_{u}|^{2p-2}\left\langle Z_{u}^{\epsilon}-Z_{u},\left(   \sigma(\hat{X}_u^{\epsilon},\sL_{\hat{X}_u^{\epsilon}})-\sigma(X_{u}^{0},\delta_{X_{u}^{0}}) \right)\dif W_{u} \right\rangle\right),\\
&&U_7:=p(2p-1)\mE\left(\sup\limits_{s\in[0,t]}\int_{0}^{s}|Z_{u}^{\epsilon}-Z_{u}|^{2p-2}\|\sigma(\hat{X}_u^{\epsilon},\sL_{\hat{X}_u^{\epsilon}})-\sigma(X_{u}^{0},\delta_{X_{u}^{0}})\|^{2}\dif u\right).
\de

Next, we estimate $U_{1}, U_{2}, U_{3}, U_{4}, U_{5}, U_{6}, U_{7}$, respectively. For $U_1, U_{3}, U_{4}, U_{5}$, by $(\bf{H}_{3})$ and $(\bf{H}_{4})$, it holds that
\be
U_{1}+U_3+U_4+U_{5}&\leq& 2pL_4\sqrt{\e}\mE\left(\int_{0}^{t}|Z_{u}^{\epsilon}-Z_{u}|^{2p-1}|Z_{u}^{\epsilon}|^2\dif u\right)\no\\
&&+2p\mE\left(\int_{0}^{t}|Z_{u}^{\epsilon}-Z_{u}|^{2p}\left\|\nabla b(X_{u}^{0},\sL_{\hat{X}_{u}^{\epsilon}})\right\|\dif u\right)\no\\
&&+2pL_4\mE\left(\int_{0}^{t}|Z_{u}^{\epsilon}-Z_{u}|^{2p-1}\mW_2(\sL_{\hat{X}_u^{\epsilon}},\delta_{X_{u}^{0}})|Z_u|\dif u\right)\no\\
&&+2p\mE\left(\int_{0}^{t}|Z_{u}^{\epsilon}-Z_{u}|^{2p-1}|D^{L}b(X_{u}^{0},\delta_{X_{u}^{0}})(X_{u}^{0})|\mE|Z_{u}^{\epsilon}-Z_{u}|\dif u\right)\no\\
&\leq&C\mE\int_{0}^{t}|Z_{u}^{\epsilon}-Z_{u}|^{2p}\dif u+C\e^p\mE\int_{0}^{t}|Z_{u}^{\epsilon}|^{4p}\dif u\no\\
&&+C\e^p\mE\int_{0}^{t}|Z_u|^{4p}\dif u,
\label{u1345es}
\ee
where the fact $\mW_{2}(\sL_{\hat{X}_{u}^{\epsilon}},\delta_{X_{u}^{0}})\leq\sqrt{\epsilon}(\mE|Z_{u}^{\epsilon}|^{2})^{\frac{1}{2}}$ is used. For $U_{2}$, $(\bf{H}_{4})$ implies that
\be
U_{2}&\leq&C\mE\int_{0}^{t}\left|\frac{b(X_u^{0},\sL_{\hat{X}_u^{\epsilon}})-b(X_u^{0},\delta_{X_u^{0}})}{\sqrt{\epsilon}}-\mE\left\langle D^{L}b(X_{u}^{0},\delta_{X_{u}^{0}})(X_{u}^{0}),Z_{u}^{\epsilon}\right\rangle\right|^{2p}\dif u\no\\
&&+C\mE\int_{0}^{t}|Z_{u}^{\epsilon}-Z_{u}|^{2p}\dif u\no\\
&=&C\mE\int_{0}^{t}\left|\int_0^1\mE\left\langle D^{L}b(X_{u}^{0},\sL_{R_{u}(r)})(R_{u}(r)),Z_{u}^{\epsilon}\right\rangle\dif r-\mE\left\langle D^{L}b(X_{u}^{0},\delta_{X_{u}^{0}})(X_{u}^{0}),Z_{u}^{\epsilon}\right\rangle\right|^{2p}\dif u\no\\
&&+C\mE\int_{0}^{t}|Z_{u}^{\epsilon}-Z_{u}|^{2p}\dif u\no\\
&\leq&C\mE\int_{0}^{t}(\mE|Z_{u}^{\epsilon}|^{2}) ^{p}\left( \int_{0}^{1}\left((\mE|R_{u}(r)-X_{u}^{0}|^{2} )^{\frac{1}{2}} +\mW_{2}(\sL_{R_{u}(r)},\delta_{X_{u}^{0}}) \right)\dif r\right)^{2p}\dif u\no\\
&&+C\mE\int_{0}^{t}|Z_{u}^{\epsilon}-Z_{u}|^{2p}\dif u\no\\
&\leq&C\e^p\int_{0}^{t}\mE|Z_{u}^{\epsilon}|^{4p}\dif u+C\mE\int_{0}^{t}|Z_{u}^{\epsilon}-Z_{u}|^{2p}\dif u,
\label{u2es}
\ee
where $R_{u}(r):=X_{u}^{0}+r(\hat{X}_{u}^{\epsilon}-X_{u}^{0})$, $r\in[0,1]$. For $U_{6}$, using the BDG inequality, the Young inequality, the H\"older inequality and $(\bf{H}_{3})$, we have that
\be
U_{6}&\leq& C\mE\left[\int_{0}^{t}|Z_{u}^{\epsilon}-Z_{u}|^{4p-2}\|   \sigma(\hat{X}_u^{\epsilon},\sL_{\hat{X}_u^{\epsilon}})-\sigma(X_{u}^{0},\delta_{X_{u}^{0}})\|^{2}\dif u\right]^{\frac{1}{2}}\no\\
&\leq&C\mE\left[\sup\limits_{u\in[0,t]}|Z_{u}^{\epsilon}-Z_{u}|^{p}\left(\int_{0}^{t}|Z_{u}^{\epsilon}-Z_{u}|^{2p-2}\|\sigma(\hat{X}_u^{\epsilon},\sL_{\hat{X}_u^{\epsilon}})-\sigma(X_{u}^{0},\delta_{X_{u}^{0}})\|^{2}\dif u\right)^{\frac{1}{2}}\right]\no\\
&\leq& \frac{1}{2}\mE\left(\sup\limits_{u\in[0,t]}|Z_{u}^{\epsilon}-Z_{u}|^{2p}\right)+C\mE\int_{0}^{t}|Z_{u}^{\epsilon}-Z_{u}|^{2p-2}\|\sigma(\hat{X}_u^{\epsilon},\sL_{\hat{X}_u^{\epsilon}})-\sigma(X_{u}^{0},\delta_{X_{u}^{0}})\|^{2}\dif u\no\\
&\leq& \frac{1}{2}\mE\left(\sup\limits_{u\in[0,t]}|Z_{u}^{\epsilon}-Z_{u}|^{2p}\right)+C\mE\int_{0}^{t}|Z_{u}^{\epsilon}-Z_{u}|^{2p}\dif u\no\\
&&+C\mE\int_{0}^{t} \left(|\hat{X}_{u}^{\epsilon}-X_{u}^{0}|+\mW_{2}(\sL_{\hat{X}_u^{\epsilon}},\delta_{X_{u}^{0}})\right)^{2p}\dif u\no\\
&\leq& \frac{1}{2}\mE\left(\sup\limits_{u\in[0,t]}|Z_{u}^{\epsilon}-Z_{u}|^{2p}\right)+ C\mE \int_{0}^{t}|Z_{u}^{\epsilon}-Z_{u}|^{2p} \dif u+C\epsilon^{p}\int_{0}^{t}\mE|Z_{u}^{\epsilon}|^{2p}\dif u.
\label{u6es}
\ee
By the similar deduction as that of (\ref{u6es}), it holds that
\be
U_{7}&\leq&C\mE\int_{0}^{t}|Z_{u}^{\epsilon}-Z_{u}|^{2p-2}\| \sigma(\hat{X}_u^{\epsilon},\sL_{\hat{X}_u^{\epsilon}})-\sigma(X_{u}^{0},\delta_{X_{u}^{0}})\|^{2}\dif u\no\\
&\leq& C\mE \int_{0}^{t}|Z_{u}^{\epsilon}-Z_{u}|^{2p} \dif u+C\epsilon^{p}\int_{0}^{t}\mE|Z_{u}^{\epsilon}|^{2p}\dif u.
\label{u7es}
\ee

Combining (\ref{u1345es})-(\ref{u7es}) with (\ref{allesti}), we get that
\ce
\mE\left(\sup\limits_{s\in[0,t]}|Z_{s}^{\epsilon}-Z_{s}|^{2p}\right)&\leq& C\epsilon^{p}+C\e^p\mE\int_{0}^{t}(|Z_{u}^{\epsilon}|^{4p}+|Z_u|^{4p})\dif u\\
&&+C\int_{0}^{t}\mE\left(\sup\limits_{s\in[0,u]}|Z_{s}^{\epsilon}-Z_{s}|^{2p}\right)\dif u.
\de
Thus, from the Gronwall inequality and Lemma \ref{auxilemm1}, it follows that
$$
\mE\left(\sup\limits_{t\in[0,T]}|Z_{t}^{\epsilon}-Z_{t}|^{2p}\right)\leq C\epsilon^{p}.
$$
The proof is complete.

\section{The moderate deviation principle for multivalued McKean-Vlasov SDEs}\label{mdp}

In this section, we study the moderate deviation principle for multivalued McKean-Vlasov SDEs.

For $\epsilon>0$, consider the following equation:
\be\left\{\begin{array}{l}
\dif \frac{\bar{X}_t^{\epsilon}- X_t^{0}}{a(\epsilon)}\in \ -A(\frac{\bar{X}_t^{\epsilon}- X_t^{0}}{a(\epsilon)})\dif t+ \frac{b(\bar{X}_t^{\epsilon},\sL_{\bar{X}_t^{\epsilon}})-b(X_t^{0},\delta_{X_t^{0}})}{a(\epsilon)}\dif t+\frac{\sqrt{\epsilon}\sigma(\bar{X}_t^{\epsilon},\sL_{\bar{X}_t^{\epsilon}})}{a(\epsilon)}\dif W_t, \quad t\in[0,T],\\
\frac{\bar{X}_0^{\epsilon}- X_0^{0}}{a(\epsilon)}=0,
\end{array}
\label{mdpequa}
\right.
\ee
where $X_{\cdot}^{0}$ solves Eq.(\ref{eq3}), i.e.
\ce\left\{\begin{array}{l}
\dif X_t^{0}\in \ -A(X_t^{0})\dif t+  b(X_t^{0},\delta_{X_t^{0}})\dif t,\\
X^0_0=\xi,
\end{array}
\right.
\de
and $a(\e)$ satisfies
\be
a(\epsilon)\rightarrow0,\quad \frac{\epsilon}{a^{2}(\epsilon)}\rightarrow0\quad as\quad \epsilon\rightarrow 0.
\label{eaecon}
\ee
So, under ($\bf{H}_{3}$), we know that Eq.(\ref{mdpequa}) has a unique solution $(\frac{\bar{X}_{\cdot}^{\epsilon}- X_{\cdot}^{0}}{a(\e)}, \bar{K}_{\cdot}^\e)$ (c.f. \cite[Theorem 3.5]{G}). Set $\bar{Y}_{t}^{\epsilon}:=\frac{\bar{X}_{t}^{\epsilon}-X_{t}^{0}}{a(\epsilon)}$, and then the moderate deviation principle for Eq.(\ref{eq1}) means that $\bar{Y}_{\cdot}^{\epsilon}$ satisfies the large deviation principle. To assure this, we make the following assumption:
\begin{enumerate}[($\bf{H}'_{3}$)]
	\item $b$ and $\sigma$ are continuous and  satisfy for $(x,\mu), (x_1,\mu_1), (x_2,\mu_2)\in\mR^{d}\times{\cP_{2}(\mR^d)}$:
	\ce
	&&\|\nabla b(x,\mu)\|\leq L'_3, \quad \|D^{L}b(x,\mu)\|_{T_{\mu,2}} \leq L'_3, \quad |b(0,\delta_{0})|\leq L'_3,\\
    &&\|\sigma(x_1,\mu_1)-\sigma(x_2,\mu_2)\|\leq L'_3(|x_1-x_2|+\mW_2(\mu_1,\mu_2)),\\
   &&\|\sigma(x,\mu)\|\leq L'_3,
    \de
    where $L'_3>0$ is a constant.
\end{enumerate}

\br
We mention that $(\bf{H}'_{3})$ is stronger than $(\bf{H}_{3})$.
\er

Note that $\sL_{\bar{X}_t^{\epsilon}}$ will converge to $\d_{X_t^0}$ as $\e\rightarrow 0$. Thus, we replace $\sL_{\bar{X}_t^{\epsilon}}$ by $\d_{X_t^0}$ and construct an approximation equation of Eq.(\ref{mdpequa}) as follows:
\be\left\{\begin{array}{l}
\dif \frac{\widetilde{X}_t^{\epsilon}- X_t^{0}}{a(\epsilon)}\in \ -A(\frac{\widetilde{X}_t^{\epsilon}- X_t^{0}}{a(\epsilon)})\dif t+ \frac{b(\widetilde{X}_t^{\epsilon},\d_{X_t^0})-b(X_t^{0},\delta_{X_t^{0}})}{a(\epsilon)}\dif t+\frac{\sqrt{\epsilon}\sigma(\tilde{X}_t^{\epsilon},\delta_{X_t^{0}})}{a(\epsilon)}\dif W_t, \quad t\in[0,T],\\
\frac{\widetilde{X}_0^{\epsilon}- X_0^{0}}{a(\epsilon)}=0.
\end{array}
\label{mdpequa1}
\right.
\ee
$(\frac{\widetilde{X}_{\cdot}^{\epsilon}- X_{\cdot}^{0}}{a(\e)}, \widetilde{K}_{\cdot}^\e)$ denotes the unique solution of the above equation (c.f. \cite[Theorem 3.5]{G}). Set $\widetilde{Y}_{t}^{\epsilon}:=\frac{\widetilde{X}_{t}^{\epsilon}-X_{t}^{0}}{a(\epsilon)}$ and then we prove that $\widetilde{Y}_{\cdot}^{\epsilon}$ satisfies the large deviation principle. 

Next, since the large deviation principle does not distinguish between exponentially equivalent families, we show the exponential equivalence of $\bar{Y}_{\cdot}^{\epsilon}$ and $\widetilde{Y}_{\cdot}^{\epsilon}$, and obtain $\bar{Y}_{\cdot}^{\epsilon}$ satisfies the large deviation principle.

\subsection{The large deviation principle for $\widetilde{Y}_{\cdot}^{\epsilon}$} 

In the subsection, we establish the large deviation principle for $\widetilde{Y}_{\cdot}^{\epsilon}$. Since an equivalent argument of the large deviation principle is the Laplace principle, we prove that $\widetilde{Y}_{\cdot}^{\epsilon}$ satisfies the Laplace principle. Here we restate the conditions for the Laplace principle.

\bco\label{cond1}
Let $\mathcal{G}^{\epsilon} : C([0,T];\mathbb{R}^m)\mapsto\mS$ be  a family of measurable mappings. There exists a
measurable mapping $\mathcal{G}^{0} : C([0,T];\mathbb{R}^m)\mapsto\mS$ such that

$(i)$ for $N\in\mathbb{N}$ and $\{h_{\e}, \e>0\}\subset\mathbf{D}_2^{N}$, $h\in \mathbf{D}_2^{N}$, if $h_{\e}\rightarrow h$ as $\e\rightarrow 0$, then
\ce
\cG^{0}\left(\int_{0}^{\cdot}h_{\e}(s)\dif s\right)\longrightarrow \cG^{0}\left(\int_{0}^{\cdot}h(s)\dif s\right).
\de

$(ii)$ for $N\in\mathbb{N}$, $\{u_{\epsilon},\epsilon>0\}\subset \mathbf{A}_{2,\e}^{N}$ and $v_\e:=\frac{u_\e}{a(\e)}\in\mathbf{A}_{2}^{N}$, if $v_{\epsilon}$ converges in distribution to $u\in\mathbf{A}_{2}^{N}$ as $\e\rightarrow 0$, then
$$
\cG^{\e}\left(W(\cdot)+\frac{1}{\sqrt{\e}}\int_{0}^{\cdot}u_{\e}(s)\dif s\right) \overset{d}{\longrightarrow} \cG^{0}\left(\int_{0}^{\cdot}u(s)\dif s\right),
$$
where $\mathbf{D}_{2,\epsilon}^{N}:=\left\{h\in\mathbb{H}:\|h\|_{\mathbb{H}}^{2}\leq Na^{2}(\epsilon) \right\}$ and $\mathbf{A}_{2,\epsilon}^{N}:=\left\{u\in\mathcal{A}: u(\omega, \cdot)\in\mathbf{D}_{2,\epsilon}^{N},a.s. \omega \right\}$.
\eco

In order to prove the Laplace principle for Eq.(\ref{mdpequa1}), we will verify Condition \ref{cond1} with
\ce
\mathbb{S}:=C([0,T],\overline{\cD(A)}), \quad \cG^{\epsilon}(W):=\widetilde{Y}_{\cdot}^{\epsilon}.
\de
Consider the following controlled multivalued McKean-Vlasov SDE:
\be\left\{\begin{array}{l}
\dif \frac{\widetilde{X}_t^{\epsilon,u}- X_t^{0}}{a(\epsilon)}\in \ -A(\frac{\widetilde{X}_t^{\epsilon,u}- X_t^{0}}{a(\epsilon)})\dif t+ \frac{b(\widetilde{X}_t^{\epsilon,u},\delta_{X_t^{0}})-b(X_t^{0},\delta_{X_t^{0}})}{a(\epsilon)}\dif t\\
\quad\quad\quad\quad\quad\quad+\frac{\sigma(\widetilde{X}_t^{\epsilon,u},\delta_{X_t^{0}})u(t)}{a(\epsilon)}
+\frac{\sqrt{\epsilon}\sigma(\widetilde{X}_t^{\epsilon,u},\delta_{X_t^{0}})}{a(\epsilon)}\dif W_t, \quad t\in[0,T], \\
\frac{\widetilde{X}_t^{\epsilon,u}- X_t^{0}}{a(\epsilon)}=0, \quad u\in\mathbf{A}_{2}^{N}.
\end{array}
\label{eq6}
\right.
\ee
By the Girsanov theorem, it holds that Eq.(\ref{eq6}) has a unique solution $(\frac{\widetilde{X}_{\cdot}^{\epsilon,u}- X_{\cdot}^{0}}{a(\epsilon)},\widetilde{K}_{\cdot}^{\epsilon,u})$. Moreover, $\widetilde{Y}^{\epsilon,u}:=\frac{\widetilde{X}_{\cdot}^{\epsilon,u}- X_{\cdot}^{0}}{a(\epsilon)}=\mathcal{G}^{\epsilon}(W(\cdot)+\frac{1}{\sqrt{\epsilon}}\int_{0}^{\cdot}u(s)\dif s)$. Let $(\widetilde{Y}_{\cdot}^{0,u},\widetilde{K}_{\cdot}^{0,u})$ solve the following multivalued McKean-Vlasov differential equation:
\be\left\{\begin{array}{l}
\dif \widetilde{Y}_{t}^{0,u}\in \ -A(\widetilde{Y}_{t}^{0,u})\dif t+ \nabla_{\widetilde{Y}_{t}^{0,u}} b(X_{t}^{0},\delta_{X_{t}^{0}})\dif t+\sigma(X_{t}^{0},\delta_{X_{t}^{0}})u(t)\dif t,\\
\widetilde{Y}_{0}^{0,u}=0, \quad u\in\mathbf{A}_{2}^{N}.
\end{array}
\label{eq5}
\right.
\ee
Thus, we define $\mathcal{G}^{0}:C([0,T];\mR^{m})\mapsto \mS$ by
$$
\mathcal{G}^{0}\left(\int_{0}^{\cdot}u(s)\dif s \right)=\widetilde{Y}^{0,u}.
$$

Next, we prove that Condition \ref{cond1} $(i)$ holds under $(\bf{H}_{3})$. 

\bl\label{mdp2}
Assume that $(\bf{H}_{3})$ holds. If $h_{\e}\rightarrow h$ in $\mathbf{D}_2^{N}$ as $\e\rightarrow0$, then $\cG^{0}(\int_{0}^{\cdot}h_{\e}(s)\dif s)$ converges to $\cG^{0}(\int_{0}^{\cdot}h(s)\dif s)$.
\el
\begin{proof}
Note that $\cG^{0}(\int_{0}^{\cdot}h_{\e}(s)\dif s)=\widetilde{Y}^{0,h_{\e}}$, $\cG^{0}(\int_{0}^{\cdot}h(s)\dif s)=\widetilde{Y}^{0,h}$ and 
\ce
&&\widetilde{Y}_{t}^{0,h_{\epsilon}}=-\widetilde{K}_{t}^{0,h_\e}+\int_{0}^{t}\nabla_{\widetilde{Y}_{s}^{0,h_{\epsilon}}} b(X_{s}^{0},\delta_{X_{s}^{0}})\dif s+\int_{0}^{t}\sigma(X_{s}^{0},\delta_{X_{s}^{0}})h_{\epsilon}(s)\dif s,\\
&&\widetilde{Y}_{t}^{0,h}=-\widetilde{K}_{t}^{0,h}+\int_{0}^{t}\nabla_{\widetilde{Y}_{s}^{0,h}} b(X_{s}^{0},\delta_{X_{s}^{0}})\dif s+\int_{0}^{t}\sigma(X_{s}^{0},\delta_{X_{s}^{0}})h(s)\dif s.
\de
Thus, by the Taylor formula, we obtain that
\ce
|\widetilde{Y}_{t}^{0,h_{\epsilon}}-\widetilde{Y}_{t}^{0,h}|^{2}&=&2\int_{0}^{t}\left\langle \widetilde{Y}_{s}^{0,h_{\epsilon}}-\widetilde{Y}_{s}^{0,h} ,\nabla_{\widetilde{Y}_{s}^{0,h_{\epsilon}}-\widetilde{Y}_{s}^{0,h}} b(X_{s}^{0},\delta_{X_{s}^{0}})\right\rangle\dif s\\
&&+2\int_{0}^{t}\left\langle   \widetilde{Y}_{s}^{0,h_{\epsilon}}-\widetilde{Y}_{s}^{0,h} ,\sigma(X_{s}^{0},\delta_{X_{s}^{0}})(h_{\epsilon}(s)-h(s))  \right\rangle \dif s\\
&&-2\int_{0}^{t}\left\langle   \widetilde{Y}_{s}^{0,h_{\epsilon}}-\widetilde{Y}_{s}^{0,h} ,\dif \widetilde{K}_{s}^{0,h_{\epsilon}}-\dif \widetilde{K}_{s}^{0,h}\right\rangle\\
&\leq&2\int_{0}^{t}\left\langle \widetilde{Y}_{s}^{0,h_{\epsilon}}-\widetilde{Y}_{s}^{0,h} ,\nabla_{\widetilde{Y}_{s}^{0,h_{\epsilon}}-\widetilde{Y}_{s}^{0,h}} b(X_{s}^{0},\delta_{X_{s}^{0}})\right\rangle\dif s\\
&&+2\int_{0}^{t}\left\langle   \widetilde{Y}_{s}^{0,h_{\epsilon}}-\widetilde{Y}_{s}^{0,h} ,\sigma(X_{s}^{0},\delta_{X_{s}^{0}})(h_{\epsilon}(s)-h(s))  \right\rangle \dif s,
\de
and
\ce
\sup\limits_{s\in[0,t]}|\widetilde{Y}_{s}^{0,h_{\epsilon}}-\widetilde{Y}_{s}^{0,h}|^{2}&\leq&2\sup\limits_{s\in[0,t]} \int_{0}^{s}\left|\left\langle \widetilde{Y}_{u}^{0,h_{\epsilon}}-\widetilde{Y}_{u}^{0,h},\nabla_{\widetilde{Y}_{u}^{0,h_{\epsilon}}-\widetilde{Y}_{u}^{0,h}} b(X_{u}^{0},\delta_{X_{u}^{0}})\right\rangle\right|\dif u\\
&&+2\sup\limits_{s\in[0,t]} \int_{0}^{s}\left|\left\langle  \widetilde{Y}_{u}^{0,h_{\epsilon}}-\widetilde{Y}_{u}^{0,h},\sigma(X_{u}^{0},\delta_{X_{u}^{0}})(h_{\epsilon}(u)-h(u))\right\rangle\right|\dif u\\
&=:&I_{1}+I_{2}.
\de

For $I_1$, from $(\bf{H}_{3})$, it follows that
\ce
I_{1}&\leq& 2\int_{0}^{t}|\widetilde{Y}_{u}^{0,h_{\epsilon}}-\widetilde{Y}_{u}^{0,h}|^2\|\nabla b(X_{u}^{0},\delta_{X_{u}^{0}})\|\dif u\\
&\leq&2L_3\int_{0}^{t}\left(\sup\limits_{s\in[0,u]} |\widetilde{Y}_{s}^{0,h_{\epsilon}}-\widetilde{Y}_{s}^{0,h}|^{2}\right)\dif u.
\de
For $I_2$, by $(\bf{H}_{3})$ and the Young inequality, one can get that
\ce
I_{2}&\leq& 2\int_{0}^{t}|\widetilde{Y}_{u}^{0,h_{\epsilon}}-\widetilde{Y}_{u}^{0,h}|\|\sigma(X_{u}^{0},\delta_{X_{u}^{0}})\||h_{\epsilon}(u)-h(u)|\dif u\\
&\leq&L_3\int_{0}^{t}|\widetilde{Y}_{u}^{0,h_{\epsilon}}-\widetilde{Y}_{u}^{0,h}|^2\dif u+L_3\int_{0}^{t}|h_{\epsilon}(u)-h(u)|^2\dif u.
\de

Finally, we obtain that
$$
\sup\limits_{s\in[0,t]}|\widetilde{Y}_{s}^{0,h_{\epsilon}}-\widetilde{Y}_{s}^{0,h}|^{2}\leq L_3\int_{0}^{T}|h_{\epsilon}(u)-h(u)|^2\dif u+3L_3\int_{0}^{t}\left(\sup\limits_{s\in[0,u]} |\widetilde{Y}_{s}^{0,h_{\epsilon}}-\widetilde{Y}_{s}^{0,h}|^{2}\right)\dif u,
$$
which together with the Gronwall inequality yields that
\ce
\sup\limits_{s\in[0,T]}|\widetilde{Y}_{s}^{0,h_{\epsilon}}-\widetilde{Y}_{s}^{0,h}|^{2}\leq L_3e^{3L_3T}\int_{0}^{T}|h_{\epsilon}(u)-h(u)|^2\dif u.
\de
The proof is complete.
\end{proof}

Now, we make preparations for justifying Condition \ref{cond1} $(ii)$.

\bl\label{auxilemm4}
Under the assumption $(\bf{H}'_{3})$, it holds that for $u_{\epsilon}\in\mathbf{A}_{2,\e}^{N}$ and $u\in\mathbf{A}_{2}^{N}$, 
\ce
\mE\left(\sup\limits_{t\in[0,T]}|\widetilde{Y}_{t}^{\epsilon,u_{\epsilon}}|^{4}\right)\leq C, \quad \mE\left(\sup\limits_{t\in[0,T]}|\widetilde{Y}_{t}^{0,u}|^{4}\right)\leq C,
\de
where $C>0$ is independent of $\epsilon$.
\el
\begin{proof}
By Eq.(\ref{eq6}), we know that
\ce
\widetilde{Y}_{t}^{\epsilon,u_{\epsilon}}&=&\int_{0}^{t}\frac{b(\widetilde{X}_s^{\epsilon,u_{\epsilon}},\delta_{X_s^{0}})-b(X_s^{0},\delta_{X_s^{0}})}{a(\epsilon)}\dif s
+\int_{0}^{t}\frac{\sigma(\widetilde{X}_s^{\epsilon,u_{\epsilon}},\delta_{X_s^{0}})u_{\epsilon}(s)}{a(\epsilon)}\dif s\\
&&+\int_{0}^{t}\frac{\sqrt{\epsilon}\sigma(\widetilde{X}_s^{\epsilon,u_{\epsilon}},\delta_{X_s^{0}})}{a(\epsilon)}\dif W_s-\widetilde{K}_{t}^{\epsilon,u_{\epsilon}}.
\de
For $\alpha\in Int(\cD(A))$, by the It\^o formula and Lemma \ref{inteineq}, it holds that
\ce
&&|\widetilde{Y}_{t}^{\epsilon,u_{\epsilon}}-\alpha|^{2}\\
&=&|\alpha|^{2}+2\int_{0}^{t}\left\langle \widetilde{Y}_{s}^{\epsilon,u_{\epsilon}}-\alpha, \frac{b(\widetilde{X}_s^{\epsilon,u_{\epsilon}},\delta_{X_s^{0}})-b(X_s^{0},\delta_{X_s^{0}})}{a(\epsilon)}  \right\rangle\dif s\no\\
&&+2\int_{0}^{t}\left\langle \widetilde{Y}_{s}^{\epsilon,u_{\epsilon}}-\alpha,  \frac{\sigma(\widetilde{X}_s^{\epsilon,u_{\epsilon}},\delta_{X_s^{0}})u_{\epsilon}(s)}{a(\epsilon)} \right\rangle\dif s\no\\
&&+2\int_{0}^{t}\left\langle \widetilde{Y}_{s}^{\epsilon,u_{\epsilon}}-\alpha, \frac{\sqrt{\epsilon}\sigma(\widetilde{X}_s^{\epsilon,u_{\epsilon}},\delta_{X_s^{0}})}{a(\epsilon)}\dif W_{s} \right\rangle\no\\
&&+\int_{0}^{t}\left\| \frac{\sqrt{\epsilon}\sigma(\widetilde{X}_s^{\epsilon,u_{\epsilon}},\delta_{X_s^{0}})}{a(\epsilon)}  \right\|^{2}\dif s-2\int_{0}^{t}\left\langle  \widetilde{Y}_{s}^{\epsilon,u_{\epsilon}}-\alpha, \dif \widetilde{K}_{s}^{\epsilon,u_{\epsilon}} \right\rangle\no\\
&\leq&|\alpha|^{2}+2\int_{0}^{t}\left\langle \widetilde{Y}_{s}^{\epsilon,u_{\epsilon}}-\alpha, \frac{b(\widetilde{X}_s^{\epsilon,u_{\epsilon}},\delta_{X_s^{0}})-b(X_s^{0},\delta_{X_s^{0}})}{a(\epsilon)}  \right\rangle\dif s\no\\
&&+2\int_{0}^{t}\left\langle \widetilde{Y}_{s}^{\epsilon,u_{\epsilon}}-\alpha,  \frac{\sigma(\widetilde{X}_s^{\epsilon,u_{\epsilon}},\delta_{X_s^{0}})u_{\epsilon}(s)}{a(\epsilon)} \right\rangle\dif s\no\\
&&+2\int_{0}^{t}\left\langle \widetilde{Y}_{s}^{\epsilon,u_{\epsilon}}-\alpha, \frac{\sqrt{\epsilon}\sigma(\widetilde{X}_s^{\epsilon,u_{\epsilon}},\delta_{X_s^{0}})}{a(\epsilon)}\dif W_{s} \right\rangle\no\\
&&+\int_{0}^{t}\left\| \frac{\sqrt{\epsilon}\sigma(\widetilde{X}_s^{\epsilon,u_{\epsilon}},\delta_{X_s^{0}})}{a(\epsilon)}  \right\|^{2}\dif s+\gamma_{2}\int_{0}^{t}|\widetilde{Y}_{s}^{\epsilon,u_{\epsilon}}-\alpha|\dif s +\gamma_{3} t-\gamma_{1}|\widetilde{K}^{\epsilon,u_{\epsilon}}  |^t_{0}\no\\
&\leq&|\alpha|^{2}+2\int_{0}^{t}\left\langle \widetilde{Y}_{s}^{\epsilon,u_{\epsilon}}-\alpha, \frac{b(\widetilde{X}_s^{\epsilon,u_{\epsilon}},\delta_{X_s^{0}})-b(X_s^{0},\delta_{X_s^{0}})}{a(\epsilon)}  \right\rangle\dif s\no\\
&&+2\int_{0}^{t}\left\langle \widetilde{Y}_{s}^{\epsilon,u_{\epsilon}}-\alpha,  \frac{\sigma(\widetilde{X}_s^{\epsilon,u_{\epsilon}},\delta_{X_s^{0}})u_{\epsilon}(s)}{a(\epsilon)} \right\rangle\dif s\no\\
&&+2\int_{0}^{t}\left\langle \widetilde{Y}_{s}^{\epsilon,u_{\epsilon}}-\alpha, \frac{\sqrt{\epsilon}\sigma(\widetilde{X}_s^{\epsilon,u_{\epsilon}},\delta_{X_s^{0}})}{a(\epsilon)}\dif W_{s} \right\rangle\no\\
&&+\int_{0}^{t}\left\| \frac{\sqrt{\epsilon}\sigma(\widetilde{X}_s^{\epsilon,u_{\epsilon}},\delta_{X_s^{0}})}{a(\epsilon)} \right\|^{2}\dif s+\gamma_{2}\int_{0}^{t}|\widetilde{Y}_{s}^{\epsilon,u_{\epsilon}}-\alpha|\dif s +\gamma_{3} t.
\de
Therefore, we get that
\ce
\mE\left(\sup\limits_{s\in[0,t]}\bigg|\widetilde{Y}_{s}^{\epsilon,u_{\epsilon}}-\alpha|^{4}\right)
&\leq& 6\left(|\alpha|^{2}+(\gamma_{2}+\gamma_{3})T\right)^2+6\gamma^2_{2}T\mE\int_{0}^{t}|\widetilde{Y}_{r}^{\epsilon,u_{\epsilon}}-\alpha|^4\dif r\\
&&+24\mE\left[\sup\limits_{s\in[0,t]}\left| \int_{0}^{s} \left\langle \widetilde{Y}_{r}^{\epsilon,u_{\epsilon}}-\alpha,\frac{b(\widetilde{X}_r^{\epsilon,u_{\epsilon}},\delta_{X_r^{0}})-b(X_r^{0},\delta_{X_r^{0}})}{a(\epsilon)}\right\rangle\dif r\right|^2 \right]\\
&&+24\mE \left[ \sup\limits_{s\in[0,t]}\bigg| \int_{0}^{s} \left\langle \widetilde{Y}_{r}^{\epsilon,u_{\epsilon}}-\alpha,\frac{\sigma(\widetilde{X}_r^{\epsilon,u_{\epsilon}},\delta_{X_r^{0}})u_{\epsilon}(r)}{a(\epsilon)}\right\rangle\dif r\bigg|^2\right]\\
&&+24\mE \left[ \sup\limits_{s\in[0,t]}\bigg| \int_{0}^{s} \left\langle \widetilde{Y}_{r}^{\epsilon,u_{\epsilon}}-\alpha,\frac{\sqrt{\epsilon}\sigma(\widetilde{X}_r^{\epsilon,u_{\epsilon}},\delta_{X_r^{0}})}{a(\epsilon)}\dif W_r \right\rangle\bigg|^2 \right]\\
&&+6\mE \left[ \sup\limits_{s\in[0,t]}\int_{0}^{s}\left\|\frac{\sqrt{\epsilon}\sigma(\widetilde{X}_r^{\epsilon,u_{\epsilon}},\delta_{X_r^{0}})}{a(\epsilon)} \right\|^{2}\dif r  \right]^2\\
&=& 6\left(|\alpha|^{2}+(\gamma_{2}+\gamma_{3})T\right)^2+6\gamma^2_{2}T\mE \int_{0}^{t}|\widetilde{Y}_{r}^{\epsilon,u_{\epsilon}}-\alpha|^4\dif r\\
&&+I_{1}+I_{2}+I_{3}+I_{4}.
\de

For $I_1$, by $(\bf{H}'_{3})$ and $\frac{\widetilde{X}_{s}^{\epsilon,u_{\epsilon}}-X_{s}^{0}}{a(\epsilon)}=\widetilde{Y}_{s}^{\epsilon,u_{\epsilon}}$, it holds that 
\ce
I_{1}&\leq& 12T\mE\int_{0}^{t}|\widetilde{Y}_{r}^{\epsilon,u_{\epsilon}}-\alpha|^4\dif r+12T\mE\int_{0}^{t}\bigg|\frac{b(\widetilde{X}_r^{\epsilon,u_{\epsilon}},\delta_{X_r^{0}})-b(X_r^{0},\delta_{X_r^{0}})}{a(\epsilon)}\bigg|^{4}\dif r\\
&\leq& 12T\mE\int_{0}^{t}|\widetilde{Y}_{r}^{\epsilon,u_{\epsilon}}-\alpha|^4\dif r+12TL_3^{'4}\mE\int_{0}^{t}\left|\frac{\widetilde{X}_r^{\epsilon,u_{\epsilon}}-X_r^{0}
}{a(\epsilon)}\right|^4\dif r\\
&\leq& (12T+96TL_3^{'4})\mE\int_{0}^{t}|\widetilde{Y}_{r}^{\epsilon,u_{\epsilon}}-\alpha|^4\dif r+96T^2L_3^{'4}|\a|^4.
\de

Noting that $u_{\epsilon}\in\mathbf{A}_{2,\epsilon}^{N}$, by the Young inequality and  $(\bf{H}'_{3})$, we get
\ce
I_{2}&\leq&24T\mE\int_{0}^{t}|\widetilde{Y}_{r}^{\epsilon,u_{\epsilon}}-\alpha|^2\left|\frac{\sigma(\widetilde{X}_r^{\epsilon,u_{\epsilon}},\delta_{X_r^{0}})u_{\epsilon}(r)}{a(\epsilon)} \right|^2\dif r\\
&\leq&24TL_3^2\mE\left(\sup\limits_{r\in[0,t]}|\widetilde{Y}_{r}^{\epsilon,u_{\epsilon}}-\alpha|^2\int_{0}^{t}\left|\frac{u_{\epsilon}(r)}{a(\epsilon)} \right|^2\dif r\right)\\
&\leq&\frac{1}{4}\mE\left(\sup\limits_{r\in[0,t]}|\widetilde{Y}_{r}^{\epsilon,u_{\epsilon}}-\alpha|^4\right)+C.
\de

For $I_3$, by the BDG inequality, the Young inequality and $(\bf{H}'_{3})$, it holds that
\ce
I_{3}&\leq&C\mE\left(\int_{0}^{t}|\widetilde{Y}_{r}^{\epsilon,u_{\epsilon}}-\alpha|^{2} \left\| \frac{\sqrt{\epsilon}\sigma(\widetilde{X}_r^{\epsilon,u_{\epsilon}},\delta_{X_r^{0}})}{a(\epsilon)} \right \|^{2}  \dif r\right) \\
&\leq&C\mE\left(\sup\limits_{r\in[0,t]}|\widetilde{Y}_{r}^{\epsilon,u_{\epsilon}}-\alpha|^2\int_{0}^{t}\left\| \frac{\sqrt{\epsilon}\sigma(\widetilde{X}_r^{\epsilon,u_{\epsilon}},\delta_{X_r^{0}})}{a(\epsilon)}  \right\|^{2}  \dif r\right) \\
&\leq&\frac{1}{4}\mE\left(\sup\limits_{r\in[0,t]}|\widetilde{Y}_{r}^{\epsilon,u_{\epsilon}}-\alpha|^4\right)+C\left(\frac{\epsilon}{a^2(\epsilon)}\right)^2.
\de
By the same deduction as that for $I_3$, one can obtain that
$$
I_{4}\leq 6L_3^{'4}T^2\left(\frac{\epsilon}{a^2(\epsilon)}\right)^2.
$$

Finally, taking the above estimates into consideration and assuming $\frac{\epsilon}{a^2(\epsilon)}\leq 1$, which is suitable in terms of (\ref{eaecon}), we get that 
\ce
\mE\left(\sup\limits_{s\in[0,t]}|\widetilde{Y}_{s}^{\epsilon,u_{\epsilon}}-\alpha|^{4}\right)
\leq C+C\int_{0}^{t}\mE\left(\sup\limits_{s\in[0,r]}|\widetilde{Y}_{s}^{\epsilon,u_{\epsilon}}-\alpha|^{4}\right) \dif r,
\de
which together with the Gronwall inequality implies the required estimate. 

As for $\widetilde{Y}^{0,u}$, it follows from Eq.(\ref{eq5}) that
\ce
\widetilde{Y}_{t}^{0,u}=\int_{0}^{t}\nabla_{\widetilde{Y}_{s}^{0,u}} b(X_{s}^{0},\delta_{X_{s}^{0}})\dif s+\int_{0}^{t}\sigma(X_{s}^{0},\delta_{X_{s}^{0}})u(s)\dif s-\widetilde{K}_{t}^{0,u}.
\de
For $\alpha\in Int(\cD(A))$, by the Taylor formula and Lemma \ref{inteineq} it holds that
\ce
&&|\widetilde{Y}_{t}^{0,u}-\alpha|^{2}\\
&=&|\alpha|^{2}-2\int_{0}^{t}\left\langle  \widetilde{Y}_{s}^{0,u} -\alpha,\dif \widetilde{K}_{s}^{0,u}\right\rangle+2\int_{0}^{t}\left\langle  \widetilde{Y}_{s}^{0,u}-\alpha  ,\nabla_{\widetilde{Y}_{s}^{0,u}} b(X_{s}^{0},\delta_{X_{s}^{0}})\right\rangle\dif s\\
&&+2\int_{0}^{t}\left\langle  \widetilde{Y}_{s}^{0,u}-\alpha,\sigma(X_{s}^{0},\delta_{X_{s}^{0}})u(s)\right\rangle\dif s\\
&\leq&|\alpha|^{2}+\gamma_{2} \int_{0}^{t}|\widetilde{Y}_{s}^{0,u}-\alpha| \dif s  +\gamma_{3} t   -\gamma_{1}|\widetilde{K}^{0,u}|_0^{t}+2\int_{0}^{t}\left\langle  \widetilde{Y}_{s}^{0,u}-\alpha,\nabla_{\widetilde{Y}_{s}^{0,u}} b(X_{s}^{0},\delta_{X_{s}^{0}})\right\rangle\dif s\\
&&+2\int_{0}^{t}\left\langle  \widetilde{Y}_{s}^{0,u} -\alpha ,\sigma(X_{s}^{0},\delta_{X_{s}^{0}})u(s)\right\rangle\dif s\\
&\leq&\left(|\alpha|^{2}+(\gamma_{2}+\gamma_{3})T\right)+\gamma_{2} \int_{0}^{t}|\widetilde{Y}_{s}^{0,u}-\alpha|^2\dif s+2\int_{0}^{t}\left\langle  \widetilde{Y}_{s}^{0,u}-\alpha  ,\nabla_{\widetilde{Y}_{s}^{0,u}} b(X_{s}^{0},\delta_{X_{s}^{0}})\right\rangle\dif s\\
&&+2\int_{0}^{t}\left\langle  \widetilde{Y}_{s}^{0,u}  -\alpha ,\sigma(X_{s}^{0},\delta_{X_{s}^{0}})u(s)\right\rangle\dif s.
\de

Noting that $u\in\mathbf{A}_{2}^{N}$, by the H\"older inequality and $(\bf{H}'_{3})$, we get that
\ce
&&\mE\left(\sup\limits_{s\in[0,t]}|\widetilde{Y}_{t}^{0,u}-\alpha|^{4}\right)\\
&\leq&4\left(|\alpha|^{2}+(\gamma_{2}+\gamma_{3})T\right)^2+4\g_2^2\mE\left[\sup\limits_{s\in[0,t]}\left(\int_{0}^{s}|\widetilde{Y}_{r}^{0,u}-\alpha|^2  \dif r\right)^{2}\right]\\
&&+16\mE\left[\sup\limits_{s\in[0,t]}\left|\int_{0}^{s}\left\langle \widetilde{Y}_{r}^{0,u}-\alpha  ,\nabla_{\widetilde{Y}_{r}^{0,u}} b(X_{r}^{0},\delta_{X_{r}^{0}})     \right\rangle\dif r\right|^{2}\right]\\
&&+16\mE\left[\sup\limits_{s\in[0,t]}\left|\int_{0}^{s}\left\langle \widetilde{Y}_{r}^{0,u} -\alpha ,\sigma(X_{r}^{0},\delta_{X_{r}^{0}})u(r)      \right\rangle\dif r\right|^{2}\right]\\
&\leq&4\left(|\alpha|^{2}+(\gamma_{2}+\gamma_{3})T\right)^2+4\g_2^2T\mE\left[\int_{0}^{t}|\widetilde{Y}_{r}^{0,u}-\alpha|^4  \dif r\right]\\
&&+16T\mE\int_{0}^{t}|\widetilde{Y}_{r}^{0,u}-\alpha|^{2} |\nabla_{\widetilde{Y}_{r}^{0,u}} b(X_{r}^{0},\delta_{X_{r}^{0}})|^{2}      \dif r\\
&&+16\mE\left(\int_{0}^{t}|\widetilde{Y}_{r}^{0,u}-\alpha|  |\sigma(X_{r}^{0},\delta_{X_{r}^{0}})u(r)| \dif r\right)^2 \\
&\leq&4\left(|\alpha|^{2}+(\gamma_{2}+\gamma_{3})T\right)^2+4\g_2^2T\mE\left[\int_{0}^{t}|\widetilde{Y}_{r}^{0,u}-\alpha|^4  \dif r\right]\\
&&+32TL_3^{'2}(1+|\a|^2)\mE\left[\int_{0}^{t}|\widetilde{Y}_{r}^{0,u}-\alpha|^4  \dif r\right]+32T^2L_3^{'2}|\a|^2\\
&&+16\mE\left(\int_{0}^{t}|\widetilde{Y}_{r}^{0,u}-\alpha|^{2}  \|\sigma(X_{r}^{0},\delta_{X_{r}^{0}})\|^{2} \dif r\right)\left(\int_{0}^{t}u^{2}(r)\dif r\right) \\
&\leq&C+C\int_{0}^{t}\mE\left(\sup\limits_{s\in[0,r]}|\widetilde{Y}_{s}^{0,u}-\alpha|^{4}\right)\dif r.
\de
By the Gronwall inequality, one can obtain that
$$
\mE\left(\sup\limits_{t\in[0,T]}|\widetilde{Y}_{t}^{0,u}-\a|^{4}\right)\leq C.
$$
The proof is complete.
\end{proof}

\bl\label{mdp1}
Assume that $(\bf{H}'_{3})$ and $(\bf{H}_{4})$ hold. Suppose that $\{u_{\epsilon}, \e>0\}\subset\mathbf{A}_{2,\e}^{N}$, $v_\e:=\frac{u_{\epsilon}}{a(\e)}\in\mathbf{A}_{2}^{N}$, $u\in\mathbf{A}_{2}^{N}$, $v_\e$ converges to $u$ almost surely as $\epsilon\rightarrow 0$. Then $\cG^{\epsilon}(W+\frac{1}{\sqrt{\epsilon}}\int_{0}^{\cdot}u_{\epsilon}(s)\dif s)\rightarrow \cG^{0}(\int_{0}^{\cdot}u(s)\dif s)$ in probability.
\el
\begin{proof}
Note that $\cG^{\epsilon}(W+\frac{1}{\sqrt{\epsilon}}\int_{0}^{\cdot}u_{\epsilon}(s)\dif s)=\widetilde{Y}_{\cdot}^{\epsilon,u_{\epsilon}}$, $\cG^{0}(\int_{0}^{\cdot}u(s)\dif s)=\widetilde{Y}_{\cdot}^{0,u}$ and 
\ce
&&\widetilde{Y}_{t}^{\epsilon,u_{\epsilon}}=-\widetilde{K}_{t}^{\epsilon,u_{\epsilon}}+\int_{0}^{t}\frac{b(\widetilde{X}_s^{\epsilon,u_{\epsilon}},\d_{X_s^0})-b(X_s^{0},\delta_{X_s^{0}})}{a(\epsilon)}\\
&&\qquad\qquad +\int_{0}^{t}\frac{\sigma(\widetilde{X}_s^{\epsilon,u_{\epsilon}},\d_{X_s^0})u_{\epsilon}(s)}{a(\epsilon)}\dif s+\int_{0}^{t}\frac{\sqrt{\epsilon}\sigma(\widetilde{X}_s^{\epsilon,u_{\epsilon}},\d_{X_s^0})}{a(\epsilon)}\dif W_{s},\\
&&\widetilde{Y}_{t}^{0,u}=-\widetilde{K}_{t}^{0,u}+\int_{0}^{t}\nabla_{\widetilde{Y}_{s}^{0,u}} b(X_{s}^{0},\delta_{X_{s}^{0}})\dif s+\int_{0}^{t}\sigma(X_{s}^{0},\delta_{X_{s}^{0}})u(s)\dif s.
\de
Thus, by the It\^o formula, it holds that
\ce
&&|\widetilde{Y}_{t}^{\epsilon,u_{\epsilon}}-\widetilde{Y}_{t}^{0,u}|^{2}\\ 
&=&    2\int_{0}^{t}\left\langle \widetilde{Y}_{s}^{\epsilon,u_{\epsilon}}-\widetilde{Y}_{s}^{0,u}, \frac{b(\widetilde{X}_s^{\epsilon,u_{\epsilon}},\delta_{X_{s}^{0}})-b(X_s^{0},\delta_{X_s^{0}})}{a(\epsilon)}-
\nabla_{\widetilde{Y}_{s}^{0,u}} b(X_{s}^{0},\delta_{X_{s}^{0}}) \right\rangle\dif s\\
&&+2\int_{0}^{t}\left\langle \widetilde{Y}_{s}^{\epsilon,u_{\epsilon}}-\widetilde{Y}_{s}^{0,u},\frac{\sigma(\widetilde{X}_s^{\epsilon,u_{\epsilon}},\delta_{X_{s}^{0}})u_{\epsilon}(s)}{a(\epsilon)}-\sigma(X_{s}^{0},\delta_{X_{s}^{0}})u(s) \right\rangle\dif s\\
&&+2\int_{0}^{t}\left\langle \widetilde{Y}_{s}^{\epsilon,u_{\epsilon}}-\widetilde{Y}_{s}^{0,u},\frac{\sqrt{\epsilon}\sigma(\widetilde{X}_s^{\epsilon,u_{\epsilon}},\delta_{X_{s}^{0}})}{a(\epsilon)}  \dif W_{s}\right\rangle\\
&&+\int_{0}^{t}\left\|  \frac{\sqrt{\epsilon}\sigma(\widetilde{X}_s^{\epsilon,u_{\epsilon}},\delta_{X_{s}^{0}})}{a(\epsilon)}  \right\|^{2}\dif s\\
&&-2\int_{0}^{t}\left\langle \widetilde{Y}_{s}^{\epsilon,u_{\epsilon}}-\widetilde{Y}_{s}^{0,u}, \dif \widetilde{K}_{s}^{\epsilon,u_{\epsilon}}-\dif \widetilde{K}_{s}^{0,u}\right\rangle\\
&\leq&2\int_{0}^{t}\left\langle \widetilde{Y}_{s}^{\epsilon,u_{\epsilon}}-\widetilde{Y}_{s}^{0,u}, \frac{b(\widetilde{X}_s^{\epsilon,u_{\epsilon}},\delta_{X_s^{0}})-b(X_s^{0},\delta_{X_s^{0}})}{a(\epsilon)}-
\nabla_{\widetilde{Y}_{s}^{0,u}} b(X_{s}^{0},\delta_{X_{s}^{0}}) \right\rangle\dif s\\
&&+2\int_{0}^{t}\left\langle \widetilde{Y}_{s}^{\epsilon,u_{\epsilon}}-\widetilde{Y}_{s}^{0,u},\frac{\sigma(\widetilde{X}_s^{\epsilon,u_{\epsilon}},\delta_{X_{s}^{0}})u_{\epsilon}(s)}{a(\epsilon)}-\sigma(X_{s}^{0},\delta_{X_{s}^{0}})u(s) \right\rangle\dif s\\
&&+2\int_{0}^{t}\left\langle \widetilde{Y}_{s}^{\epsilon,u_{\epsilon}}-\widetilde{Y}_{s}^{0,u},\frac{\sqrt{\epsilon}\sigma(\widetilde{X}_s^{\epsilon,u_{\epsilon}},\delta_{X_{s}^{0}})}{a(\epsilon)}  \dif W_{s}\right\rangle\\
&&+\int_{0}^{t}\left\|  \frac{\sqrt{\epsilon}\sigma(\widetilde{X}_s^{\epsilon,u_{\epsilon}},\delta_{X_{s}^{0}})}{a(\epsilon)} \right\|^{2}\dif s.
\de
Therefore, we get that
\ce
&&\mE\left(\sup\limits_{s\in[0,t]}|\widetilde{Y}_{s}^{\epsilon,u_{\epsilon}}-\widetilde{Y}_{s}^{0,u}|^{2}\right)\\
&\leq&2\mE\left[\sup\limits_{s\in[0,t]} \bigg|  \int_{0}^{s}\left\langle \widetilde{Y}_{r}^{\epsilon,u_{\epsilon}}-\widetilde{Y}_{r}^{0,u}, \frac{b(\widetilde{X}_r^{\epsilon,u_{\epsilon}},\delta_{X_r^{0}})-b(X_r^{0},\delta_{X_r^{0}})}{a(\epsilon)}-
\nabla_{\widetilde{Y}_{r}^{0,u}} b(X_{r}^{0},\delta_{X_{r}^{0}}) \right\rangle\dif r     \bigg|     \right]\\
&&+2\mE\left[\sup\limits_{s\in[0,t]}  \bigg|  \int_{0}^{s}\left\langle \widetilde{Y}_{r}^{\epsilon,u_{\epsilon}}-\widetilde{Y}_{r}^{0,u},\frac{\sigma(\widetilde{X}_r^{\epsilon,u_{\epsilon}},\delta_{X_r^{0}})u_{\epsilon}(r)}{a(\epsilon)}-\sigma(X_{r}^{0},\delta_{X_{r}^{0}})u(r) \right\rangle\dif r    \bigg|    \right] \\
&&+2\mE\left[\sup\limits_{s\in[0,t]}  \bigg|   \int_{0}^{s}\left\langle \widetilde{Y}_{r}^{\epsilon,u_{\epsilon}}-\widetilde{Y}_{r}^{0,u},\frac{\sqrt{\epsilon}\sigma(\widetilde{X}_r^{\epsilon,u_{\epsilon}},\delta_{X_r^{0}})}{a(\epsilon)}  \dif W_{r}\right\rangle    \bigg|    \right] \\
&&+\mE\left[\sup\limits_{s\in[0,t]} \int_{0}^{s}\left\|  \frac{\sqrt{\epsilon}\sigma(\widetilde{X}_r^{\epsilon,u_{\epsilon}},\delta_{X_r^{0}})}{a(\epsilon)}  \right\|^{2}\dif r \right]\\
&=:&I_{1}+I_{2}+I_{3}+I_{4}.
\de

For $I_{1}$, by the Young inequality, $(\bf{H}'_{3})$ and $(\bf{H}_{4})$, it holds that
\be
I_1&\leq&\mE\int_{0}^{t}|\widetilde{Y}_{r}^{\epsilon,u_{\epsilon}}-\widetilde{Y}_{r}^{0,u}|^2\dif r\no\\
&&+\mE\int_{0}^{t}\bigg| \frac{b(\widetilde{X}_r^{\epsilon,u_{\epsilon}},\delta_{X_r^{0}})-b(X_r^{0},\delta_{X_r^{0}})}{a(\epsilon)}-
\nabla_{\widetilde{Y}_{r}^{0,u}} b(X_{r}^{0},\delta_{X_{r}^{0}})\bigg|^{2}\dif r\no\\
&\leq&\mE\int_{0}^{t}|\widetilde{Y}_{r}^{\epsilon,u_{\epsilon}}-\widetilde{Y}_{r}^{0,u}|^2\dif r\no\\
&&+\mE\int_{0}^{t}\left| \int_{0}^{1}  \nabla_{\widetilde{Y}_{r}^{\epsilon,u_{\epsilon}}} b(X_{r}^{0}+\eta(\widetilde{X}_r^{\epsilon,u_{\epsilon}} -X_{r}^{0}),\delta_{X_{r}^{0}}) \dif \eta-\nabla_{\widetilde{Y}_{r}^{0,u}} b(X_{r}^{0},\delta_{X_{r}^{0}}) \right|^{2}\dif r\no\\
&\leq& \mE\int_{0}^{t}|\widetilde{Y}_{r}^{\epsilon,u_{\epsilon}}-\widetilde{Y}_{r}^{0,u}|^2\dif r+2\mE\int_{0}^{t} \left| \int_{0}^{1} \nabla_{\widetilde{Y}_{r}^{\epsilon,u_{\epsilon}}-\widetilde{Y}_{r}^{0,u}} b(X_{u}^{0}+\eta(\widetilde{X}_r^{\epsilon,u_{\epsilon}} -X_{r}^{0}),\delta_{X_{r}^{0}})\dif \eta \right|^{2}\dif r\no\\
&&+2\mE\int_{0}^{t}\left|\int_{0}^{1} \nabla_{\widetilde{Y}_{r}^{0,u}} b(X_{r}^{0}+\eta(\widetilde{X}_r^{\epsilon,u_{\epsilon}} -X_{r}^{0}),\delta_{X_{r}^{0}})\dif \eta  -\nabla_{\widetilde{Y}_{r}^{0,u}} b(X_{r}^{0},\delta_{X_{r}^{0}}) \right|^{2}\dif r\no\\
&\leq& (1+2L_3^{'2})\mE\int_{0}^{t}|\widetilde{Y}_{r}^{\epsilon,u_{\epsilon}}-\widetilde{Y}_{r}^{0,u}|^2\dif r+2L_4^2a^2(\e)\left(\mE\int_{0}^{T}|\widetilde{Y}_{r}^{\epsilon,u_{\epsilon}}|^4\dif r\right)^{1/2}\left(\mE\int_{0}^{T}|\widetilde{Y}_{r}^{0,u}|^4\dif r\right)^{1/2}.\no\\
\label{i1ales}
\ee

For $I_{2}$, by the Young inequality, one can obtain that
\ce
I_{2}&\leq&2\mE\left[\int_{0}^{t}|\widetilde{Y}_{r}^{\epsilon,u_{\epsilon}}-\widetilde{Y}_{r}^{0,u}|\left|\frac{\sigma(\widetilde{X}_r^{\epsilon,u_{\epsilon}},\delta_{X_{r}^{0}})u_{\epsilon}(r)}{a(\epsilon)}-\sigma(X_{r}^{0},\delta_{X_{r}^{0}})u(r) \right|\dif r \right] \\
&\leq&\mE\int_{0}^{t}|\widetilde{Y}_{r}^{\epsilon,u_{\epsilon}}-\widetilde{Y}_{r}^{0,u}|^2\dif r+2\mE\int_{0}^{t}\left|\frac{\left(\sigma(\widetilde{X}_r^{\epsilon,u_{\epsilon}},\delta_{X_{r}^{0}})-\sigma(X_{r}^{0},\delta_{X_{r}^{0}})\right)u_{\epsilon}(r) }{a(\epsilon)}\right|^2\dif r\\
&&+2\mE\int_{0}^{t}\left|\sigma(X_{r}^{0},\delta_{X_{r}^{0}})\left(\frac{u_{\epsilon}(r)}{a(\epsilon)}-u(r)\right)\right|^2\dif r\\
&=:&\mE\int_{0}^{t}|\widetilde{Y}_{r}^{\epsilon,u_{\epsilon}}-\widetilde{Y}_{r}^{0,u}|^2\dif r+I_{21}+I_{22}.
\de
Noting that $u_{\epsilon}\in A_{2,\epsilon}^{N}$, we get
\ce
I_{21}&\leq&2L^{'2}_3\mE\int_{0}^{t}|\widetilde{X}_r^{\epsilon,u_{\epsilon}}-X_{r}^{0}|^2\left|\frac{u_{\epsilon}(r)}{a(\epsilon)}\right|^2\dif r\\
&\leq&2L^{'2}_3a^2(\epsilon)N\mE\left(\sup\limits_{r\in[0,T]}|\widetilde{Y}_{r}^{\epsilon,u_{\epsilon}}|^2\right).
\de
For $I_{22}$, by $(\bf{H}'_{3})$, it holds that
\ce
I_{22}\leq 2L^{'2}_3\mE\int_{0}^{T}\left|\frac{u_{\epsilon}(r)}{a(\epsilon)}-u(r)\right|^2\dif r.
\de
Thus, we know that
\be
I_{2}&\leq& \mE\int_{0}^{t}|\widetilde{Y}_{r}^{\epsilon,u_{\epsilon}}-\widetilde{Y}_{r}^{0,u}|^2\dif r+2L^{'2}_3a^2(\epsilon)N\mE\left(\sup\limits_{r\in[0,T]}|\widetilde{Y}_{r}^{\epsilon,u_{\epsilon}}|^2\right)\no\\
&&+2L^{'2}_3\mE\int_{0}^{T}\left|\frac{u_{\epsilon}(r)}{a(\epsilon)}-u(r)\right|^2\dif r.
\label{i2ales}
\ee

By the BDG inequality and the Young inequality, it holds that
\be
I_{3}&\leq&C\mE\left(\int_{0}^{t}|\widetilde{Y}_{r}^{\epsilon,u_{\epsilon}}-\widetilde{Y}_{r}^{0,u}|^{2}\left\|\frac{\sqrt{\epsilon}\sigma(\widetilde{X}_r^{\epsilon,u_{\epsilon}},\delta_{X_{r}^{0}})}{a(\epsilon)}  \right \|^{2}\dif r\right)^{1/2}\no\\
&\leq&CL'_3\mE\left[ \frac{\sqrt{\epsilon}}{a(\epsilon)} \left(\int_{0}^{t} |\widetilde{Y}_{r}^{\epsilon,u_{\epsilon}}-\widetilde{Y}_{r}^{0,u}|^{2} \dif r\right)^{1/2}  \right]\no\\
&\leq& C\frac{\epsilon}{a^{2}(\epsilon)}+C\mE\int_{0}^{t} |\widetilde{Y}_{r}^{\epsilon,u_{\epsilon}}-\widetilde{Y}_{r}^{0,u}|^{2} \dif r.
\label{i3ales}
\ee

For $I_4$, by $(\bf{H}'_{3})$, we get
\be
I_{4}\leq L^{'2}_3T\frac{\epsilon}{a^{2}(\epsilon)}.
\label{i4ales}
\ee

Finally, combining (\ref{i1ales})-(\ref{i4ales}), one can have that
\ce
\mE\left(\sup\limits_{s\in[0,t]}|\widetilde{Y}_{s}^{\epsilon,u_{\epsilon}}-\widetilde{Y}_{s}^{0,u}|^{2}\right)&\leq& C\int_{0}^{t}\mE\left(\sup\limits_{s\in[0,r]}|\widetilde{Y}_{s}^{\epsilon,u_{\epsilon}}-\widetilde{Y}_{s}^{0,u}|^2\right)\dif r\\
&&+2L_4^2a^2(\e)\left(\mE\int_{0}^{T}|\widetilde{Y}_{r}^{\epsilon,u_{\epsilon}}|^4\dif r\right)^{1/2}\left(\mE\int_{0}^{T}|\widetilde{Y}_{r}^{0,u}|^4\dif r\right)^{1/2}\\
&&+2a^2(\e)L^{'2}_3N\mE\left(\sup\limits_{r\in[0,T]}|\widetilde{Y}_{r}^{\epsilon,u_{\epsilon}}|^2\right)\\
&&+2L^{'2}_3\mE\int_{0}^{T}\left|\frac{u_{\epsilon}(r)}{a(\epsilon)}-u(r)\right|^2\dif r+C\frac{\epsilon}{a^{2}(\epsilon)},
\de
which together with the Gronwall inequality and Lemma \ref{auxilemm4} implies that
\ce
\mE\left(\sup\limits_{s\in[0,T]}|\widetilde{Y}_{s}^{\epsilon,u_{\epsilon}}-\widetilde{Y}_{s}^{0,u}|^{2}\right)\leq \left(Ca^2(\e)+2L^{'2}_3\mE\int_{0}^{T}\left|\frac{u_{\epsilon}(r)}{a(\epsilon)}-u(r)\right|^2\dif r+C\frac{\epsilon}{a^{2}(\epsilon)}\right)e^{CT}.
\de
By the dominated convergence theorem, we know that $\mE\int_{0}^{T}\left|\frac{u_{\epsilon}(r)}{a(\epsilon)}-u(r)\right|^2\dif r\rightarrow 0$ as $\e\rightarrow 0$. Noting that $\frac{\epsilon}{a^{2}(\epsilon)}\rightarrow 0$ and $a^{2}(\epsilon)\rightarrow 0$ as $\epsilon\rightarrow 0$, we conclude that
\ce
\mE\left(\sup\limits_{s\in[0,T]}|\widetilde{Y}_{s}^{\epsilon,u_{\epsilon}}-\widetilde{Y}_{s}^{0,u}|^{2}\right)\rightarrow 0, ~\mbox{as}~ \epsilon\rightarrow 0.
\de
The proof is complete.
\end{proof}

\bt\label{th2}
Assume that $(\bf{H}'_{3})$ and $(\bf{H}_{4})$ hold. Then the family $\{\widetilde{Y}^{\epsilon}\}$ satisfies the Laplace principle in $\mS:=C([0,T],\overline{\cD(A)})$ with the rate function given by
\ce
I(x)=\frac{1}{2}\inf\limits_{h\in\mH, x=\widetilde{Y}^{0,h}}\|h\|_{\mH}^{2}.
\de
\et
\begin{proof}
By Lemma \ref{mdp2}, we know that Condition \ref{cond1} $(i)$ holds. 

Next, we verify Condition \ref{cond1} $(ii)$. For $\epsilon\in(0,1)$ and $\{u_{\epsilon}\}\subset\mathbf{A}_{2,\e}^{N}$, $u\in\mathbf{A}_{2}^{N}$, let $\frac{u_{\epsilon}}{a(\e)}$ converge to $u$ in distribution. By the Skorohod theorem, there exists a 
probability space $(\check{\Omega}, \check{\sF}, \check{\mP})$, and ${\bf D}_2^N$-valued random variables $\{\frac{\check{u}_{\epsilon}}{a(\e)}\}$, $\check{u}$ and a $m$-dimensional Brownian motion $\check{W}$ defined on it such that

(i) $\sL_{(\frac{\check{u}_{\epsilon}}{a(\e)},\check{W})}=\sL_{(\frac{u_{\epsilon}}{a(\e)},W)}$ and $\sL_{\check{u}}=\sL_{u}$;

(ii) $\frac{\check{u}_{\epsilon}}{a(\e)}$ converges to $\check{u}$ almost surely.

Besides, by Lemma \ref{mdp1} we have that
\ce
\mathcal{G}^{\epsilon}\left(\check{W}(\cdot)+\frac{1}{\sqrt{\epsilon}}\int_{0}^{\cdot}\check{u}_{\epsilon}(s)\dif s\right)\overset{\check{\mP}}\longrightarrow \mathcal{G}^{0}\left(\int_{0}^{\cdot}\check{u}(s)\dif s \right).
\de
Thus, it holds that
\ce
\mathcal{G}^{\epsilon}\left(\check{W}(\cdot)+\frac{1}{\sqrt{\epsilon}}\int_{0}^{\cdot}\check{u}_{\epsilon}(s)\dif s\right)\overset{d}\longrightarrow \mathcal{G}^{0}\left(\int_{0}^{\cdot}\check{u}(s)\dif s \right),
\de
which yields that
\ce
\mathcal{G}^{\epsilon}\left(W(\cdot)+\frac{1}{\sqrt{\epsilon}}\int_{0}^{\cdot}u_{\epsilon}(s)\dif s\right)\overset{d}\longrightarrow \mathcal{G}^{0}\left(\int_{0}^{\cdot}u(s)\dif s \right).
\de
Therefore, Condition \ref{cond1} $(ii)$ is right.

Lastly, by Theorem \ref{ldpbase}, we know that the family $\{\widetilde{Y}^{\epsilon}\}$ satisfies the Laplace principle. The proof is complete.
\end{proof}

\subsection{The large deviation principle for $\bar{Y}_{\cdot}^{\epsilon}$}

In the subsection, we prove that $\bar{Y}_{\cdot}^{\epsilon}$ and $\widetilde{Y}_{\cdot}^{\epsilon}$ are exponentially equivalent and then obtain the large deviation estimate for $\bar{Y}_{\cdot}^{\epsilon}$.

\bl\label{ee}
Assume that $(\bf{H}'_{3})$ and $(\bf{H}_{4})$ hold. Then it holds that for any $\delta>0$,
\be
\limsup\limits_{\epsilon\rightarrow0}\epsilon\log\left(\mP\left\{\sup\limits_{t\in[0,T]}|\bar{Y}_{t}^{\epsilon}-\widetilde{Y}_{t}^{\epsilon}|\geq\delta \right\} \right)=-\infty.
\label{ee0}
\ee
\el
\begin{proof}
Note that $\bar{Y}^{\epsilon}=\frac{\bar{X}_{t}^{\epsilon}-X_{t}^{0}}{a(\epsilon)}, \widetilde{Y}^{\epsilon}=\frac{\widetilde{X}_{t}^{\epsilon}-X_{t}^{0}}{a(\epsilon)}$ satisfy the following equations, respectively,
\ce
&&\bar{Y}_{t}^{\epsilon}=\int_0^t\frac{b(\bar{X}_s^{\epsilon},\sL_{\bar{X}_s^{\epsilon}})-b(X_s^{0},\delta_{X_s^{0}})}{a(\epsilon)}\dif s+\int_0^t\frac{\sqrt{\epsilon}\sigma(\bar{X}_s^{\epsilon},\sL_{\bar{X}_s^{\epsilon}})}{a(\epsilon)}\dif W_s-\bar{K}^\e_t,\\
&&\widetilde{Y}_{t}^{\epsilon}=\int_0^t\frac{b(\widetilde{X}_s^{\epsilon},\d_{X_s^0})-b(X_s^{0},\delta_{X_s^{0}})}{a(\epsilon)}\dif s+\int_0^t\frac{\sqrt{\epsilon}\sigma(\tilde{X}_s^{\epsilon},\delta_{X_s^{0}})}{a(\epsilon)}\dif W_s-\tilde{K}^\e_t.
\de
Thus, set $G_{t}:=\bar{Y}_{t}^{\epsilon}-\widetilde{Y}_{t}^{\epsilon}$, and then $G_{\cdot}$ satisfies the following equation
\ce
G_{t}=\int_{0}^{t}b_{s}\dif s+\sqrt{\epsilon}\int_{0}^{t}\sigma_{s}\dif W_s-\bar{K}^\e_t+\tilde{K}^\e_t,
\de
where 
$$
b_{s}:=\frac{b(\bar{X}_s^{\epsilon},\sL_{\bar{X}_s^{\epsilon}})-b(\widetilde{X}_s^{\epsilon},\d_{X_s^0})}{a(\epsilon)}, \quad \sigma_{s}:=\frac{\sigma(\bar{X}_s^{\epsilon},\sL_{\bar{X}_s^{\epsilon}})-\sigma(\tilde{X}_s^{\epsilon},\delta_{X_s^{0}})}{a(\epsilon)}.
$$
Moreover, by $(\bf{H}'_{3})$, it holds that 
\ce
|b_{s}|&\leq&L'_{3}\left(\frac{|\bar{X}_s^{\epsilon}-\widetilde{X}_s^{\epsilon}|}{a(\epsilon)}+\frac{\mW_{2}(\sL_{\bar{X}_s^{\epsilon}},\d_{X_s^0})}{a(\epsilon)}\right)\leq \sqrt{2}L'_{3}\left(|G_{s}|^{2}+\rho^{2}(\epsilon)\right)^{\frac{1}{2}},\\
|\sigma_{s}|&\leq&L'_{3}\left(\frac{|\bar{X}_s^{\epsilon}-\widetilde{X}_s^{\epsilon}|}{a(\epsilon)}+\frac{\mW_{2}(\sL_{\bar{X}_s^{\epsilon}},\d_{X_s^0})}{a(\epsilon)}\right)\leq \sqrt{2}L'_{3}\left(|G_{s}|^{2}+\rho^{2}(\epsilon)\right)^{\frac{1}{2}},
\de
where $\rho^{2}(\epsilon)=\sup\limits_{t\in[0,T]}\mE|\bar{Y}_t^{\epsilon}|^{2}$.

Next, we choose a $R>0$ such that $|X_{t}^{0}|<R+1$ for $t\in[0,T]$. Then define a stopping time as follows:
\ce
\tau_{R}=\inf\{t\geq 0: |\bar{Y}_{t}^{\epsilon}|\vee|\widetilde{Y}_{t}^{\epsilon}|\geq R+1\}\wedge T.
\de
Thus, by the similar deduction to that of \cite[Lemma 5.6.18]{D}, it holds that for any $\delta>0$ and $0<\epsilon\leq1$, 
\ce
\epsilon \log\left(\mP\left\{\sup\limits_{t\in[0,\tau_{R}]}|G_{t}|\geq\delta \right\} \right)\leq C+\log\left(\frac{\rho^{2}(\epsilon)}{\rho^{2}(\epsilon)+\delta ^{2}} \right).
\de
By the same deduction to that of Lemma \ref{xex0}, one can prove that $\lim\limits_{\epsilon\rightarrow0}\rho^{2}(\epsilon)=0$, which implies that 
\be
\lim\limits_{\epsilon\rightarrow0}\epsilon \log\left(\mP\left\{\sup\limits_{t\in[0,\tau_{R}]}|G_{t}|\geq\delta \right\} \right)=-\infty.
\label{ee1}
\ee

Besides, we investigate $\{\tau_{R}\leq T\}$. Set $\eta_{R}:=\{t\geq 0: |\widetilde{Y}_{t}^{\epsilon}|\geq R\}$, and then it holds that
\ce
\mP\{ \tau_{R}\leq T\}&\leq&\mP\{|\bar{Y}_{\tau_{R}}^{\epsilon}|=R+1 \}+\mP\{|\widetilde{Y}_{\tau_{R}}^{\epsilon}|=R+1\}\\
&\leq&\mP\left\{\sup\limits_{t\in[0,\tau_{R}]}|G_{t}|\geq\frac{1}{2},|\bar{Y}_{\tau_{R}}^{\epsilon}|=R+1 \right\}+\mP\left\{\sup\limits_{t\in[0,\tau_{R}]}|G_{t}|<\frac{1}{2},|\bar{Y}_{\tau_{R}}^{\epsilon}|=R+1 \right\}\\
&&+\mP\{\eta_{R}\leq T\}\\
&\leq&\mP\left\{\sup\limits_{t\in[0,\tau_{R}]}|G_{t}|\geq\frac{1}{2} \right\}+2\mP\{\eta_{R}\leq T\}.
\de
On the one hand, from (\ref{ee1}), it follows that 
$$
\lim\limits_{\epsilon\rightarrow0}\epsilon \log\left(\mP\left\{\sup\limits_{t\in[0,\tau_{R}]}|G_{t}|\geq\frac{1}{2} \right\}\right)=-\infty. 
$$
On the other hand, note that
\ce
\limsup\limits_{\epsilon\rightarrow0}\epsilon \log\mP\{\eta_{R}\leq T  \} &\leq&\limsup\limits_{\epsilon\rightarrow0}\epsilon \log\bigg(\mP\bigg\{\sup\limits_{t\in[0,T]}|\widetilde{Y}_{t}^{\epsilon}|\geq R\bigg\}  \bigg)\\
&\leq&-\inf_{\stackrel{h\in\mH, x=\widetilde{Y}^{0,h},}{\sup\limits_{t\in[0,T]}|x_t|\geq R}}\frac{1}{2}\|h\|_{\mH}^{2},
\de
where the last inequality is based on Theorem \ref{th2}. Now we observe the last term. By Eq.(\ref{eq5}), it holds that
\ce
\widetilde{Y}_{t}^{0,h}=-\widetilde{K}_{t}^{0,h}+\int_{0}^{t}\nabla_{\widetilde{Y}_{s}^{0,h}} b(X_{s}^{0},\delta_{X_{s}^{0}})\dif s+\int_{0}^{t}\sigma(X_{s}^{0},\delta_{X_{s}^{0}})h(s)\dif s.
\de
Using the similar method to one in the proof of Lemma \ref{auxilemm4}, we obtain that
\ce
\sup\limits_{t\in[0,T]}|\widetilde{Y}_{t}^{0,h}|^{2}\leq 2|\a|^2+2\left(C+L'_3\int_0^T|h(t)|^2\dif t\right)e^{CT}<\infty,
\de
where $\alpha\in Int(\cD(A))$. From this, it follows that 
$$
\left\{h\in\mH, x=\widetilde{Y}^{0,h},\sup\limits_{t\in[0,T]}|x_t|\geq R\right\}\longrightarrow\emptyset, ~\mbox{as}~ R\rightarrow \infty,
$$  
which implies that
$$
-\lim\limits_{R\rightarrow\infty}\inf_{\stackrel{h\in\mH, x=\widetilde{Y}^{0,h},}{\sup\limits_{t\in[0,T]}|x_t|\geq R}}\frac{1}{2}\|h\|_{\mH}^{2}=-\infty,
$$
and furthermore
\ce
\lim\limits_{R\rightarrow\infty}\limsup\limits_{\epsilon\rightarrow0}\epsilon \log\bigg(\mP\{\eta_{R}\leq T  \}  \bigg)=-\infty.
\de
Taking the above estimates into consideration, we know that
\be
\lim\limits_{R\rightarrow\infty}\limsup\limits_{\epsilon\rightarrow0}\epsilon \log\bigg(\mP\{\tau_{R}\leq T  \}  \bigg)=-\infty.
\label{ee2}
\ee

Finally, note that
$$
\left\{\sup\limits_{t\in[0,T]}|G_t|\geq\delta\right\}\subset\left\{\tau_{R}\leq T  \right\}\bigcup\left\{\sup\limits_{t\in[0,\tau_{R}]}|G_t|\geq\delta \right\}.
$$
Combining (\ref{ee1}) and (\ref{ee2}), we get
\ce
\limsup\limits_{\epsilon\rightarrow0}\epsilon \log\mP\left(\sup\limits_{t\in[0,T]}|Y_{t}^{\epsilon}-\widetilde{Y}_{t}^{\epsilon}|\geq\delta  \right)=-\infty.
\de
The proof is over.
\end{proof}

Now, by Theorem \ref{th2} and Lemma \ref{ee}, we draw the following conclusion which is the main result in the section.

\bt\label{th3}
Assume that $(\bf{H}'_{3})$ and $(\bf{H}_{4})$ hold. Then the family $\{Y^{\epsilon}\}$ satisfies the Laplace principle in $\mS:=C([0,T],\overline{\cD(A)})$ with the rate function given by
\ce
I(x)=\frac{1}{2}\inf\limits_{h\in\mH, x=\widetilde{Y}^{0,h}}\|h\|_{\mH}^{2}.
\de
\et

\bigskip

\textbf{Acknowledgements:}

Two authors would like to thank Doctor Wei Hong for pointing out a mistake in the origin version.


\begin{thebibliography}{999}

\bibitem{cepa} E. C\'epa: \'Equations diff\'erentielles stochasticques multivoques, {\it Lect Notes in Math S\'eminaire Prob XXIX}, Springer, Berlin, 1995, 86-107.
	
\bibitem{cepaa} E. C\'epa: Probleme de Skorohod Multivoque, {\it Ann. Prob.}, 26(1998)500-532.

\bibitem{CHI} H. Chi: Multivalued stochastic McKean-Vlasov equation, {\it Acta Math. Sci.}, 34B(2014)1731-1740.

\bibitem{D} A. Dembo and O. Zeitouni:  {\it Large Deviations Techniques and Applications}, vol. 38, Springer, Berlin, 2010.

\bibitem{DQ} X. Ding and H. Qiao: Stability for stochastic McKean-Vlasov equations with non-Lipschitz coefficients, {\it SIAM J. Control Optim.}, 59(2021)887-905.

\bibitem{DQ2} X. Ding and H. Qiao: Euler-Maruyama approximations for stochastic McKean-Vlasov equations with non-Lipschitz coefficients, {\it Journal of Theoretical Probability}, 34(2021)1408-1425.

\bibitem{rst} G. Dos Reis, W. Salkekd and J. Tugaut: Freidlin-Wentzell LDP in path space for McKean-Vlasov equations and the functional iterated logarithm law, {\it Ann. Appl. Probab.}, 29(2019)1487-1540.

\bibitem{de} P. Dupuis and R. S. Ellis: {\it A Weak Convergence Approach to the Theory of Large Deviations}, Wiley, New York, 1997.

\bibitem{G} J. Gong and H. Qiao: The stability for multivalued McKean-Vlasov SDEs with non-Lipschitz coefficients, https://arxiv.org/abs/2106.12080.

\bibitem{Kac} M. Kac: {\it Foundations of kinetic theory}. In Proceedings of the Third Berkeley Symposium on Mathematical Statistics and Probability, 1954-1955, vol. III, pages 171-197. University of California Press, Berkeley and Los Angeles, 1956.

\bibitem{lwyz} Y. Li, R. Wang, N. Yao, and S. Zhang: A moderate deviation principle for stochastic Volterra equation, {\it Stat. Probab. Lett.}, 122 (2017)79-85.

\bibitem{lq} M. Liu and H. Qiao: Parameter estimation of path-dependent McKean-Vlasov stochastic differential equations, {\it Acta Mathematica Scientia}, 42B(2022)876-886.

\bibitem{lszz} W. Liu, Y. Song, J. Zhai and T. Zhang: Large and moderate deviation principles for McKean-Vlasov SDEs with jumps, https://arxiv.org/abs/2011.08403.

\bibitem{V} V. Maroulas: Large deviations of infinite dimensional stochastic systems with jumps, {\it Mathematika}, 57(2011)175-192.

\bibitem{m} T. M. Manil: Moderate deviation principle for the 2D stochastic convective Brinkman-Forchheimer equations, {\it Stochastics}, 93(2021)1052-1106.

\bibitem{RW} Y. Ren and J. Wang: Large deviation for mean-field stochastic differential equations with subdifferential operator, {\it Stochastic Analysis and applications}, 34(2016)318-338.

\bibitem{Re} P. Ren and F.-Y. Wang: Bismut formula for lions derivative of distribution dependent SDEs and applications, {\it J. Differ. Equ} 267(2019)4745-4777.

\bibitem{Ren} J. Ren and J. Wu: The optimal control problem associated with multi-valued stochastic differential equations with jumps, {\it Nonlinear Anal},  86(2013)30-51.

\bibitem{rwz} J. Ren, J. Wu and H. Zhang, General large deviations and functional iterated logarithm law for multivalued stochastic differential equations, {\it J. Theoret. Probab.} 28(2015)550-586.

\bibitem{RWZ1} J. Ren, J. Wu and X. Zhang: Exponential ergodicity of non-Lipschitz multivalued stochastic differential equations, {\it Bull. Sci. Math}, 134(2010)391-404.

\bibitem{RXZ} J. Ren, S. Xu and X. Zhang: Large deviations for multivalued stochastic differential equations, {\it J. Theoret. Probab}, 23(2010)1142-1156.

\bibitem{Suo} Y. Suo and C. Yuan: Central limit theorem and moderate deviation principle for McKean-Vlasov SDEs, {\it Acta Applicandae Mathematicae}, (2021)175:16.

\bibitem{w} L. Wu: Moderate deviations of dependent random variables related to CLT, {\it Ann. Probab.}, 23(1995)420-445.

\bibitem{zh} H. Zhang: Moderate deviation principle for multivalued stochastic differential equations, {\it Stochastic and Dynamics}, 20(2020)1-30.

\bibitem{ZXCH} X. Zhang: Skorohod problem and multivalued stochastic evolution equations in Banach spaces, {\it Bull Sci Math}, 131(2007)175-217.
\end{thebibliography}
\end{document}